\documentclass[5p,sort&compress]{elsarticle}
\usepackage{subfigure,graphicx}
\usepackage{subcaption}
\usepackage{amsmath}
\usepackage{yhmath}
\usepackage{amsfonts}
\usepackage{amssymb}
\usepackage{epsf}
\usepackage{color}
\usepackage{mathtools}
\usepackage{placeins}
\usepackage{booktabs}
\usepackage{wrapfig}
\usepackage{caption}
\usepackage{psfrag}
\usepackage{enumitem}
\usepackage{adjustbox}
\usepackage{epstopdf}
\usepackage{dcolumn}
\usepackage{amsmath}
\usepackage{amsfonts}
\usepackage{amssymb}
\usepackage{colortbl}
\usepackage{color}
\usepackage{bbm}
\usepackage{float}
\usepackage{amsthm}
\usepackage{tabu}
\usepackage{placeins}
\usepackage{siunitx}
\usepackage{hyperref}

\theoremstyle{remark}
\newtheorem{remark}{Remark}
\biboptions{sort&compress}

\def\bfb{{\bf b}}

\def\bfx{{\bf x}}
\def\bfy{{\bf y}}

\def\bfA{{\bf A}}

\def\bfB{{\bf B}}
\def\bfC{{\bf C}}

\def\bfI{{\bf I}}

\def\bfN{{\bf N}}

\def\bfS{{\bf S}}
\def\bfT{{\bf T}}
\def\bfX{{\bf X}}
\def\bfF{{\bf F}}

\def\bfD{{\bf D}}
\def\bfe{{\bf e}}

\def\bfG{\textbf{G}}

\def\bfw{\textbf{w}}

\def\Atan{\mbox{\boldmath$\mathcal{A}$}}

\def\l{\lambda}

\def\e0{\varepsilon_0}

\def\oI{\overline{I}}

\def\obfB{\overline{\bfB}}
\def\obfC{\overline{\bfC}}

\begin{document}

\begin{frontmatter}

\title{Abaqus implementation of a large family of finite viscoelasticity models}

\author[NU]{Victor Lef\`evre}
\ead{victor.lefevre@northwestern.edu}

\author[Illinois]{Fabio Sozio}
\ead{fsozio@illinois.edu}

\author[Illinois]{Oscar Lopez-Pamies}
\ead{pamies@illinois.edu}

\address[NU]{Department of Mechanical Engineering, Northwestern University, Evanston, IL 60208, USA}

\address[Illinois]{Department of Civil and Environmental Engineering, University of Illinois, Urbana--Champaign, IL 61801, USA}

\begin{abstract}

In this paper, we introduce an Abaqus UMAT subroutine for a family of constitutive models for the viscoelastic response of isotropic elastomers of any compressibility --- including fully incompressible elastomers --- undergoing finite deformations. The models can be chosen to account for a wide range of non-Gaussian elasticities, as well as for a wide range of nonlinear viscosities. From a mathematical point of view, the structure of the models is such that the viscous dissipation is characterized by an internal variable $\bfC^v$, subject to the physically-based constraint $\det\bfC^v=1$, that is solution of a nonlinear first-order ODE in time. This ODE is solved by means of an explicit Runge-Kutta scheme of high order capable of preserving the constraint $\det\bfC^v=1$ identically. The accuracy and convergence of the code is demonstrated numerically by comparison with an exact solution for several of the Abaqus built-in hybrid finite elements, including the simplicial elements C3D4H and C3D10H and the hexahedral elements C3D8H and C3D20H. The last part of this paper is devoted to showcasing the capabilities of the code by deploying it to compute the homogenized response of a bicontinuous rubber blend.

\vspace{0.2cm}

\keyword{Elastomers; Rubber blends; Finite deformations; Hybrid finite elements; Stable ODE solvers}
\endkeyword

\end{abstract}

\end{frontmatter}

\section{Introduction}

Over the past three decades, primarily because of their superior numerical tractability, constitutive models based on internal variables \cite{Sidoroff74,LeTallec93,RG98,BB98,Linder11,KLP16,Ravi20,Cohen21} have established themselves as the preferred choice over models based on hereditary integrals \cite{GreenRivlin57,PipkinRogers68,Lockett72} to describe the mechanical \emph{dissipative} response of polymers, hydrogels, soft biological tissues, and other soft organic materials.

In this context, the objective of this work is to put forth an Abaqus UMAT (user material) subroutine for a family of internal-variable-based constitutive models for the finite viscoelastic response of elastomers, that introduced by Kumar and Lopez-Pamies \cite{KLP16}. Such models can be derived within the so-called two-potential framework \cite{HN75,Ziegler87,KLP16} and hence are characterized by two thermodynamic potentials: $i$) a free-energy function $\psi$ that serves to characterize how the material stores energy through elastic deformation and $ii$) a dissipation potential $\phi$ that serves to characterize how the material dissipates energy through viscous deformation. A distinguishing advantage of this approach is that it allows to enforce material frame indifference, material symmetry, and entropy imbalance from the outset in a straightforward manner \cite{KLP16,Yavari24}.

The focus of this paper is on models for isotropic elastomers of any compressibility, including fully incompressible elastomers. By construction, the models allow to describe the elasticity of these materials in terms of an initial bulk modulus $\kappa\in(0,+\infty)$ and two non-Gaussian stored-energy functions $\Psi^{{\rm Eq}}$ and $\Psi^{{\rm NEq}}$ of choice. They also allow to describe their viscosity  in terms of a nonlinear viscosity function $\eta$ of choice.

From a mathematical point of view, the structure of the models is such that they provide the first Piola-Kirchhoff stress $\bfS(\bfX,t)$ at any material point $\bfX$ and time $t$ explicitly in terms of the deformation gradient $\bfF(\bfX,t)$, a pressure field $q(\bfX,t)$, and an internal variable $\bfC^v(\bfX,t)$, subject to the constraint\footnote{Physically, the constraint $\det\bfC^v=1$ describes that viscous dissipation in elastomers is an isochoric process.} $\det\bfC^v(\bfX,t)=1$, that is solution of a nonlinear first-order ordinary differential equation (ODE) in time. The fact that the internal variable $\bfC^v$ satisfies the non-convex constraint $\det\bfC^v=1$ poses one of the main difficulties. Indeed, as is well-known from the classical literature on finite plasticity \cite{Simo92}, commonly
used time integration schemes are unable to deliver solutions that satisfy such a constraint. In this work, we handle this challenge by making use of a time integration scheme based on an explicit fifth-order-accurate Runge-Kutta integrator capable of preserving the constraint $\det\bfC^v=1$ identically \cite{Lawson66,KLP16,WLPM23}.

The organization of the paper is as follows. We begin in Section \ref{Sec: The problem} by formulating the finite viscoelastostatics problem of interest in this work. The family of finite viscoelasticity models under consideration are introduced in Subsection \ref{Sec: Constitutive}, while the final set of the governing equations that they lead to for the deformation field $\bfy(\bfX,t)$, the pressure field $q(\bfX,t)$, and the internal variable $\bfC^v(\bfX,t)$ is presented in strong form in Subsection \ref{Sec: Gov Eq}. In Section \ref{Sec: Discretizations}, we present the discretization in time and space of the weak form of the governing equations. In particular, we make use of a FD (finite difference) discretization of time and a FE (finite element) discretization of space. In Section \ref{Sec: UMAT}, we describe the inputs needed for an Abaqus UMAT subroutine in order to solve in Abaqus the discretized governing equations laid out in Section \ref{Sec: Discretizations}. In Section \ref{Sec: Convergence}, we demonstrate the accuracy and convergence of the proposed Abaqus UMAT subroutine. Finally, in Section \ref{Sec: Performance}, we showcase the capabilities of the subroutine by solving a problem of both fundamental and practical interest, that of the homogenization of the finite viscoelastic response of a bicontinuous rubber blend.

\section{The problem}\label{Sec: The problem}

\subsection{Initial configuration and kinematics}\label{Sec: Kinematics}

Consider a body made of an elastomer that in its initial configuration, at time $t=0$, occupies the open domain $\mathrm{\Omega}_0\subset \mathbb{R}^3$, with boundary $\partial\mathrm{\Omega}_0$ and outward unit normal $\bfN$. We identify material points by their initial position vector $\bfX\in\mathrm{\Omega}_0$. At a later time $t\in(0,T]$, in response to the boundary conditions and body forces described in Subsection \ref{Sec: BCs} below, the position vector $\bfX$ of a material point occupies a new position $\bfx$ specified by an invertible mapping $\bfy$ from $\mathrm{\Omega}_0$ to the current configuration $\mathrm{\Omega}(t)\subset \mathbb{R}^3$. We write
\begin{equation*}
\bfx=\bfy(\bfX, t)
\end{equation*}
and the associated deformation gradient and Lagrangian velocity fields at $\bfX\in\Omega_0$ and $t\in(0,T]$ as
\begin{equation*}
\bfF(\bfX, t)=\nabla\bfy(\bfX,t)=\frac{\partial \bfy}{\partial \bfX}(\bfX,t)
\end{equation*}
and
\begin{equation*}
\textbf{V}(\bfX,t)=\dot{\bfy}(\bfX, t)= \frac{\partial \bfy}{\partial t}(\bfX,t).
\end{equation*}
Throughout, we shall use the ``dot'' notation to denote the material time derivative (i.e., with $\bfX$ held fixed) of field quantities.

\subsection{Constitutive behavior}\label{Sec: Constitutive}

As anticipated in the Introduction, the focus of this work is on viscoelastic isotropic elastomers whose isothermal mechanical behavior is described in terms of two thermodynamic potentials, a free energy function of the form \cite{KLP16}
\begin{equation}
\psi(\bfF,\bfC^v)=\psi^{{\rm Eq}}(\bfF)+\psi^{{\rm NEq}}(\bfF,\bfC^v)\label{psi_body}
\end{equation}
that describes how the elastomer stores energy through elastic deformation and a dissipation potential of the form

\begin{equation}
\phi(\bfF,\bfC^v,\dot{\bfC}^v)=\left\{\begin{array}{ll}
\hspace{-0.1cm}\dfrac{1}{2}\dot{\bfC}^v\cdot\Atan(\bfF,\bfC^v)\dot{\bfC}^v& {\rm if}\, \det\bfC^v=1\\ \\
\hspace{-0.1cm}+\infty & {\rm otherwise}\end{array}\right. \label{phi_body}
\end{equation}
that describes how the elastomer dissipates energy through viscous deformation. In these expressions, the symmetric second-order tensor $\bfC^v$ is an internal variable of state that stands for a measure of the ``viscous part'' of the deformation gradient $\bfF$, $\psi^{{\rm Eq}}$ is a non-negative function that characterizes the elastic energy storage in the elastomer at states of thermodynamic equilibrium, the non-negative function $\psi^{{\rm NEq}}$ characterizes the additional elastic energy storage at non-equilibrium states (i.e., the part of the energy that gets dissipated eventually), and the fourth-order tensor $\Atan$ characterizes the (deviatoric) viscosity of the elastomer; see Fig. \ref{Fig1} below for a rheological representation of $\psi^{{\rm Eq}}$, $\psi^{{\rm NEq}}$, and $\phi$.

For consistency with the type of constitutive relations that can be implemented as Abaqus UMAT subroutines, we consider equilibrium and non-equilibrium stored-energy functions of the form
\begin{equation}\label{Potentials-I1-J}
\left\{\begin{array}{l}
\psi^{{\rm Eq}}(\bfF)=\Psi^{{\rm Eq}}(\oI_1)+\dfrac{\kappa}{2}(J-1)^2\vspace{0.4cm}\\
\psi^{{\rm NEq}}(\bfF,\bfC^v)=\Psi^{{\rm NEq}}(\oI^e_1)\end{array}\right.
\end{equation}
and viscosity tensors of the form
\begin{equation}\label{Viscosity-A}
\mathcal{A}_{ijkl}(\bfF,\bfC^v)=\dfrac{\eta(\oI_1^e,\oI_2^e,I_1^v)}{2}{C^{v}}^{-1}_{il}{C^{v}}^{-1}_{jk}.
\end{equation}
In these expressions,
\begin{align*}
\left\{\begin{array}{l}
I_1= {\rm tr}\,\bfC \vspace{0.2cm}\\
%I_2=\dfrac{1}{2}\left[I_1^2-{\rm tr}\,\bfC^2\right] \vspace{0.2cm}\\
J=\sqrt{\det\bfC}\vspace{0.2cm}\\
\oI_1={\rm tr}\,\obfC=J^{-2/3} I_1 \vspace{0.2cm}\\
%\oI_2=\dfrac{1}{2}\left[\oI_1^2-{\rm tr}\,\obfC^2\right]=J^{-4/3} I_2 \vspace{0.2cm}\\
I_1^v={\rm tr}\,\bfC^v\vspace{0.2cm}\\
%I_2^v=\dfrac{1}{2}\left[(I_1^v)^2-{\rm tr}\,(\bfC^v)^2\right]\vspace{0.2cm}\\
I_1^e={\rm tr}\left(\bfC{\bfC^{v}}^{-1}\right)\vspace{0.2cm}\\
I_2^e=\dfrac{1}{2}\left[\left(\bfC\cdot{\bfC^v}^{-1}\right)^2-{\bfC^v}^{-1}\bfC\cdot\bfC{\bfC^v}^{-1}\right]\vspace{0.2cm}\\
\oI_1^e=J^{-2/3} I^e_1\vspace{0.2cm}\\
\oI_2^e=J^{-4/3} I^e_2\end{array}\right. ,
\end{align*}
where $\bfC=\bfF^T\bfF$, $\obfC=J^{-2/3}\bfC$, $\bfC^v={\bfF^v}^T\bfF^v$, and $\Psi^{{\rm Eq}}$, $\Psi^{{\rm NEq}}$, $\eta$ are non-negative material functions of their arguments, while $\kappa$ is a non-negative material constant. The latter describes the compressibility of the elastomer in a monotonically increasing fashion, in particular, the larger the value of $\kappa$ the more incompressible the elastomer is. The choice $\kappa=+\infty$ corresponds to a fully incompressible elastomer. In this formulation, the material constant $\kappa$ is in fact the initial bulk modulus of the elastomer in the limit of small deformations as $\bfF\rightarrow\bfI$.

Granted the two thermodynamic potentials (\ref{psi_body}) and (\ref{phi_body}), it follows that the first Piola-Kirchhoff stress tensor $\bfS$ at any material point $\bfX\in\Omega_0$ and time $t\in[0,T]$ is expediently given by the relation \cite{KLP16}
\begin{equation}
\bfS(\bfX,t)=\frac{\partial \psi}{\partial\bfF}(\bfF,\bfC^v),\label{S-gen}
\end{equation}
where $\bfC^v$ is implicitly defined by the evolution equation
\begin{equation}
\left\{\begin{array}{l}\dfrac{\partial \psi}{\partial \bfC^v}(\bfF,\bfC^v)+\dfrac{\partial \phi}{\partial \dot{\bfC}^v}(\bfF,\bfC^v,\dot{\bfC}^v)={\bf0}\\ \\
\bfC^v(\bfX,0)=\bfI\end{array}\right. .\label{Evolution-gen}
\end{equation}
Making explicit use of the representations (\ref{Potentials-I1-J})-(\ref{Viscosity-A}), the constitutive relation (\ref{S-gen})-(\ref{Evolution-gen}) specializes to
\begin{multline}
\bfS(\bfX,t)=2J^{-2/3}\Psi^{{\rm Eq}}_{\oI_1}\bfF+2J^{-2/3}\Psi^{{\rm NEq}}_{\oI^e_1}\bfF{\bfC^v}^{-1}-\\
\dfrac{2}{3}J^{-2/3}\left(I_1\Psi^{{\rm Eq}}_{\oI_1}+I^e_1\Psi^{{\rm NEq}}_{\oI^e_1}\right)\bfF^{-T}+\kappa(J-1)J\bfF^{-T},\label{S-I1-J}
\end{multline}
where $\bfC^v$ is defined implicitly as the solution of the evolution equation\footnote{Here, it is worth remarking that the stress-deformation relation (\ref{S-I1-J}) depends on $\bfC^v$ only through its inverse ${\bfC^v}^{-1}$ and that, in turn, the ODE (\ref{Evolution-I1-J})$_1$ can be rewritten solely in terms of ${\bfC^v}^{-1}$: $\dot{\wideparen{ {\bfC^v}^{-1} }} = -{\bf{G}} ( {\bfC^{v-1}}\bfC{\bfC^v}^{-1} , {\bfC^v}^{-1} )$. When implementing (\ref{S-I1-J})-(\ref{Evolution-I1-J}) numerically, one can hence write the equations in terms of $\bfD^v\equiv{\bfC^v}^{-1}$ and thereby circumvent having to perform inversions of the internal variable. This is in fact how we have coded (\ref{S-I1-J})-(\ref{Evolution-I1-J}) in the proposed UMAT subroutine.}
\begin{equation}
\left\{\begin{array}{l}\dot{\bfC}^v(\bfX,t)=\bfG(\bfC,\bfC^v)\vspace{0.2cm}\\
\hspace{1.45cm}\equiv\dfrac{2J^{-2/3}\Psi^{{\rm NEq}}_{\oI^e_1}}{\eta}\left[\bfC-\dfrac{1}{3}\left(\bfC\cdot{\bfC^v}^{-1}\right)\bfC^v\right] \\ \\
\bfC^v(\bfX,0)=\bfI\end{array}\right. .  \label{Evolution-I1-J}
\end{equation}
In these last expressions and below, we make use of the notation $f_x={\rm d} f(x)/{\rm d} x$ and $f_{xx}={\rm d}^2 f(x)/{\rm d} x^2$.

The following remarks are in order:
\begin{remark}
The Cauchy stress tensor $\bfT(\bfx,t)=J^{-1}\bfS\bfF^T$ at any spatial point $\bfx\in\Omega(t)$ and time $t\in[0,T]$ is given by
\begin{multline*}
\bfT(\bfx,t)=2J^{-5/3}\Psi^{{\rm Eq}}_{\oI_1}\bfF\bfF^T+2J^{-5/3}\Psi^{{\rm NEq}}_{\oI^e_1}\bfF{\bfC^v}^{-1}\bfF^T-\\
\left(\dfrac{2}{3}J^{-5/3}I_1\Psi^{{\rm Eq}}_{\oI_1}+\dfrac{2}{3}J^{-5/3}I^e_1\Psi^{{\rm NEq}}_{\oI^e_1}-\kappa(J-1)\right)\bfI.
\end{multline*}
When decomposed into deviatoric and volumetric contributions, it reads
\begin{equation*}
\bfT(\bfx,t)={\rm dev}\,\bfT+p\bfI
\end{equation*}
with
\begin{equation}\label{Tvol}
p=\dfrac{1}{3}{\rm tr}\, \bfT=\kappa(J-1)
\end{equation}
and
\begin{multline*}
{\rm dev}\,\bfT=\underbrace{2J^{-5/3}\Psi^{{\rm Eq}}_{\oI_1}\bfF\bfF^T-\dfrac{2}{3}J^{-5/3}I_1\Psi^{{\rm Eq}}_{\oI_1}\bfI}_{{\rm dev}\,\bfT^{{\rm Eq}}}+\\
\underbrace{2J^{-5/3}\Psi^{{\rm NEq}}_{\oI^e_1}\bfF{\bfC^v}^{-1}\bfF^T-\dfrac{2}{3}J^{-5/3}I^e_1\Psi^{{\rm NEq}}_{\oI^e_1}\bfI}_{{\rm dev}\,\bfT^{{\rm NEq}}}.
\end{multline*}
\end{remark}
\begin{remark}
The stored-energy functions $\Psi^{{\rm Eq}}$ and $\Psi^{{\rm NEq}}$ in the constitutive relation (\ref{S-I1-J})-(\ref{Evolution-I1-J}) are arbitrary functions of the invariants $\oI_1$ and $\oI^e_1$. They can thus be chosen as any of the numerous non-Gaussian $I_1$-based models available in the literature  \cite{Yeoh93,arruda1993three,Gent96,beatty2003average,Benam21,LP10}. In the representative results that we present in Sections \ref{Sec: Convergence} and \ref{Sec: Performance} below, we make use of the two-term Lopez-Pamies stored-energy functions \cite{LP10}
\begin{equation}\label{WLP-Eq}
\Psi^{{\rm Eq}}(\oI_1)=\displaystyle\sum_{r=1}^2\dfrac{3^{1-\alpha_r}}{2\alpha_r}\mu_r\left[\,\oI_1^{\alpha_r}-3^{\alpha_r}\right]
\end{equation}
and
\begin{equation}\label{WLP-NEq}
\Psi^{{\rm NEq}}(\oI^e_1)=\displaystyle\sum_{r=1}^2\dfrac{3^{1-a_r}}{2a_r}m_r\left[\,{\oI^e_1}^{a_r}-3^{a_r}\right],
\end{equation}
where $\mu_r$, $\alpha_r$, $m_r$, $a_r$ $(r=1,2)$ are material constants.
\end{remark}
\begin{remark}
The viscosity function $\eta$ in the evolution equation (\ref{Evolution-I1-J}) is an arbitrary isotropic function of the invariants
$\oI^e_1$, $\oI^e_2$, $\oI^v_1$. It can thus be chosen as needed to describe the viscosity of the elastomer at hand. Typically, elastomers exhibit a strongly nonlinear viscosity of shear-thinning type; see, e.g., \cite{Gent62,Lion06,KLP16,Khayat18,GSKLP21,SGLP23}.  In the representative results that we present in Sections \ref{Sec: Convergence} and \ref{Sec: Performance} below, we make use of the shear-thinning viscosity function introduced by Kumar and Lopez-Pamies \cite{KLP16}:
\begin{equation}\label{eta-KLP}
\eta(\oI_1^e,\oI_2^e,I_1^v)=\eta_{\infty}+\dfrac{\eta_0-\eta_{\infty}+K_1\left[{I_1^v}^{\beta_1}-3^{\beta_1}\right]}{1+\left(K_2 J_2^{{\rm NEq}}\right)^{\beta_2}}
\end{equation}
with
\begin{align*}
J_2^{{\rm NEq}}=&\dfrac{1}{2}{\rm dev}\,\bfT^{{\rm NEq}}\cdot{\rm dev}\,\bfT^{{\rm NEq}}\\
=&4J^{-2}\left(\dfrac{(\oI_1^e)^2}{3}-\oI_2^e\right)\left(\Psi^{{\rm NEq}}_{\oI^e_1}\right)^2,
\end{align*}
where $\eta_0$, $\eta_{\infty}$, $K_1$, $K_2$, $\beta_1$, $\beta_2$ are material constants.
\end{remark}
\begin{remark}
The ratio
\begin{align}\label{tau}
\tau=\dfrac{\eta}{2\Psi^{{\rm NEq}}_{\oI^e_1}}
\end{align}
in the evolution equation (\ref{Evolution-I1-J}) describes the spectrum of material time scales characteristic of the elastomer at hand.
\end{remark}
\begin{remark}
The constitutive relation (\ref{S-I1-J})-(\ref{Evolution-I1-J}) is nothing more than a generalization of the classical Zener or standard solid model \cite{Zener48} to the setting of finite deformations \cite{KLP16}. Figure \ref{Fig1} illustrates its rheological representation.
\begin{figure}[H]
   \centering \includegraphics[width=2.7in]{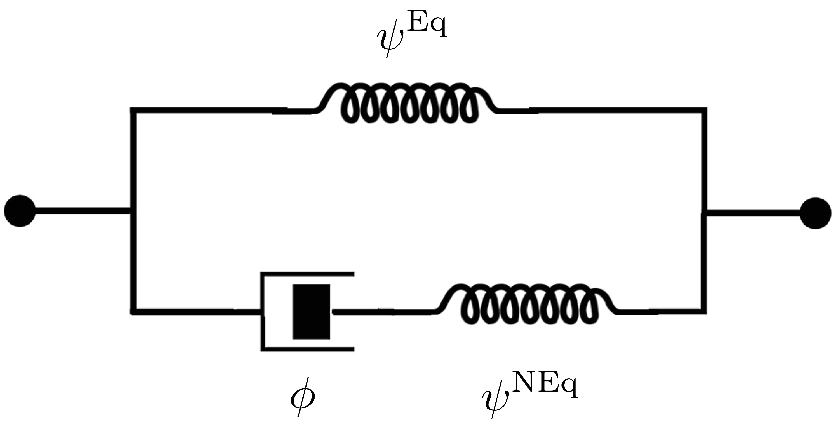}
   \caption{\small Rheological representation of the two-potential model (\ref{S-I1-J})-(\ref{Evolution-I1-J}).}
   \label{Fig1}
\end{figure}
\end{remark}
\begin{remark}
The constitutive relation (\ref{S-I1-J})-(\ref{Evolution-I1-J}) includes two fundamental models as limiting cases. The first one, which corresponds to setting the viscosity function either to $\eta=0$ or $\eta=+\infty$, is that of a \emph{non-Gaussian elastic solid}. The second one, which corresponds to setting the equilibrium and non-equilibrium energies to $\psi^{{\rm Eq}}=0$ and $\psi^{{\rm NEq}}=+\infty$, is that of an \emph{incompressible non-Newtonian fluid}.

To see the specialization to the non-Gaussian elastic solid, note that when $\eta=0$, the solution to the evolution equation (\ref{Evolution-I1-J}) is simply $\bfC^v=\bfC$ and the first Piola-Kirchhoff stress tensor (\ref{S-I1-J}) reduces to the form $\bfS(\bfX,t)=2J^{-2/3}\Psi^{{\rm Eq}}_{\oI_1}\bfF-p_1J\bfF^{-T}$. Similarly, when $\eta\rightarrow+\infty$, the solution to the evolution equation (\ref{Evolution-I1-J}) is $\bfC^v=\bfI+O(\eta^{-1})$ and the first Piola-Kirchhoff stress tensor (\ref{S-I1-J}) reduces to $\bfS(\bfX,t)=2J^{-2/3}\Psi^{{\rm Eq}}_{\oI_1}\bfF+2J^{-2/3}\Psi^{{\rm NEq}}_{\oI^e_1}\bfF-p_2J\bfF^{-T}$.

On the other hand, to see the specialization to the incompressible non-Newtonian fluid, note that when $\Psi^{{\rm Eq}}=0$ and $\Psi^{{\rm NEq}}=m[\oI_1^e-3]/2$ with $m\rightarrow+\infty$, the solution to the evolution equation (\ref{Evolution-I1-J}) is given $\bfC^v=\bfC+m^{-1}(-\eta\dot{\bfC}+q_1\bfC)+O(m^{-2})$ and the first Piola-Kirchhoff stress tensor (\ref{S-I1-J}) reduces to $\bfS(\bfX,t)=\eta(\dot{\bfF}\bfF^{-1}\bfF^{-T}+\bfF^{-T}\dot{\bfF}^T\bfF^{-T})-q_2\bfF^{-T}$; in these last two expressions, $q_1$ and $q_2$ are arbitrary hydrostatic pressures associated with the incompressibility constraint of the material. Accordingly, the Cauchy stress tensor $\bfT=\bfS\bfF^{T}$ specializes to $\bfT(\bfx,t)=2\eta\bfD-q_2\bfI$, where $\bfD=1/2(\dot{\bfF}\bfF^{-1}+\bfF^{-T}\dot{\bfF}^T)$ is the rate of deformation tensor.
\end{remark}

\subsection{Boundary conditions and body forces}\label{Sec: BCs}

The external stimuli applied to the body comprise both prescribed mechanical boundary data and body forces in the bulk. Specifically, on a portion $\partial\Omega_0^{\mathcal{D}}$ of the boundary $\partial\Omega_0$ the deformation field $\bfy(\bfX,t)$ is taken to be given by a known
function $\overline{\bfy}(\bfX,t)$, while the complementary part of the boundary $\partial\Omega_0^{\mathcal{N}}=\partial\Omega_0\setminus \partial\Omega_0^{\mathcal{D}}$
is subjected to a prescribed nominal traction $\overline{\textbf{s}}(\bfX,t)$. That is,
\begin{equation*}
\left\{\begin{array}{ll}
\bfy(\bfX,t)=\overline{\bfy}(\bfX,t), & (\bfX,t)\in\partial\Omega_0^{\mathcal{D}}\times [0,T]\\ \\
\bfS(\bfX,t)\bfN=\overline{\textbf{s}}(\bfX,t), & (\bfX,t)\in\partial\Omega_0^{\mathcal{N}}\times [0,T]\end{array}\right. . %\label{BCs}
\end{equation*}
Throughout $\Omega_0$, we also consider that the body is subjected to a mechanical body force with density
\begin{equation*}
\bfb(\bfX,t),\qquad (\bfX,t)\in\Omega_0\times [0,T] . %\label{Source}
\end{equation*}

\subsection{Governing equations: The deformation-based formulation}\label{Sec: Gov Eq 0}

Absent inertia, the relevant equations of balance of linear and angular momenta read simply as ${\rm Div}\,\bfS+\bfb=\textbf{0}$ and $\bfS\bfF^T=\bfF\bfS^T$ for $(\bfX,t)\in\Omega_0\times [0,T]$. The latter is automatically satisfied by virtue of the objectivity of the functions (\ref{Potentials-I1-J}) and so the governing equations for the response of the body reduce to the following initial-boundary-value problem:
\begin{equation}
\left\{\begin{array}{l}{\rm Div}\left[2J^{-2/3}\Psi^{{\rm Eq}}_{\oI_1}\nabla\bfy+2J^{-2/3}\Psi^{{\rm NEq}}_{\oI^e_1}\nabla\bfy{\bfC^v}^{-1}-\right.\vspace{0.2cm}\\
\hspace{0.6cm}\dfrac{2}{3}J^{-2/3}\left(I_1\Psi^{{\rm Eq}}_{\oI_1}+I^e_1\Psi^{{\rm NEq}}_{\oI^e_1}\right)\nabla\bfy^{-T}+\vspace{0.2cm}\\
\hspace{0.6cm}\left.\kappa(J-1)J\nabla\bfy^{-T}\right]+\bfb={\bf0}, \qquad (\bfX,t)\in\Omega_0\times[0,T]\\ \\
\bfy(\bfX,t)=\overline{\bfy}(\bfX,t), \qquad (\bfX,t)\in\partial\Omega^{\mathcal{D}}_0\times[0,T] \\ \\
\left[2J^{-2/3}\Psi^{{\rm Eq}}_{\oI_1}\nabla\bfy+2J^{-2/3}\Psi^{{\rm NEq}}_{\oI^e_1}\nabla\bfy{\bfC^v}^{-1}-\right.\vspace{0.2cm}\\
\dfrac{2}{3}J^{-2/3}\left(I_1\Psi^{{\rm Eq}}_{\oI_1}+I^e_1\Psi^{{\rm NEq}}_{\oI^e_1}\right)\nabla\bfy^{-T}+\vspace{0.2cm}\\
\left.\kappa(J-1)J\nabla\bfy^{-T}\right]\bfN=\overline{\textbf{s}}(\bfX,t), \quad (\bfX,t)\in\partial\Omega^{\mathcal{N}}_0\times[0,T] \\ \\
\bfy(\bfX,0)=\bfX, \qquad \bfX\in\Omega_0
\end{array}\right. \label{Equilibrium-PDE}
\end{equation}
coupled with the evolution equation
\begin{equation}
\left\{\begin{array}{l}
\hspace{-0.15cm}\dot{\bfC^v}(\bfX,t)=\bfG\left(\nabla\bfy^T\nabla\bfy,\bfC^v\right), \quad (\bfX,t)\in\Omega_0\times[0,T]\\ \\
\hspace{-0.15cm}\bfC^v(\bfX,0)=\bfI,\qquad\bfX\in\Omega_0\end{array}\right. \label{Evolution-ODE}
\end{equation}
for the deformation field $\bfy(\bfX,t)$ and the internal variable $\bfC^v(\bfX,t)$.

\subsection{Governing equations: The hybrid formulation}\label{Sec: Gov Eq}

In order to deal with nearly or fully incompressible elastomers, it is convenient to deal \emph{not} with equations (\ref{Equilibrium-PDE})-(\ref{Evolution-ODE}) directly, but with an alternative set of governing equations wherein a pressure field (and not just the deformation field) is also an unknown. Much like in the simpler setting of elasticity \cite{CTLPP15,CTLPP16,LLP17}, as outlined next, the derivation of such a mixed set of equations hinges on the introduction of an appropriate Legendre transform \cite{GSKLP21}.

Consistent with the way in which Abaqus deals with nearly incompressible materials, consider the following function
\begin{equation}\label{W-Leg}
W(\oI_1,\oI_1^e,J)=\Psi^{{\rm Eq}}(\oI_1)+\Psi^{{\rm NEq}}(\oI^e_1)+\dfrac{\kappa}{2}(J-1)^2
\end{equation}
alongside its partial Legendre transform
\begin{equation}
W^\star(\oI_1,\oI_1^e,q)=\max_{J}\left\{q(J-1)-W(\oI_1,\oI_1^e,J)\right\}.\label{Legendre-W}
\end{equation}
Since $W(\oI_1,\oI_1^e,J)$ is convex in its third argument, it readily follows that
\begin{align*}
W(\oI_1,\oI_1^e,J)&=\left(W^\star\right)^{\star}(\oI_1,\oI_1^e,J)\\
                     &=\max_{q}\left\{q(J-1)-W^\star(\oI_1,\oI_1^e,q)\right\}.
%\label{Legendre-twice}
\end{align*}
In turn, it follows that the first Piola-Kirchhoff stress tensor (\ref{S-gen}) can be rewritten in terms of the dual function (\ref{Legendre-W}) as
\begin{align}
\left\{\begin{array}{l}\bfS(\bfX,t)=-\dfrac{\partial W^\star}{\partial\bfF}(\oI_1,\oI_1^e,q)+q J\bfF^{-T}\\ \\
{\rm with}\\ \\
J=1+\dfrac{\partial W^\star}{\partial q}(\oI_1,\oI_1^e,q)\end{array}\right.
\label{Legendre-1piola-v0}
\end{align}

Making direct use of the specific form (\ref{W-Leg}) of the function $W$, a straightforward calculation shows that its partial Legendre transform (\ref{Legendre-W}) is given by
\begin{equation*}
W^\star(\oI_1,\oI_1^e,q)=\dfrac{q^2}{2\kappa}-\Psi^{{\rm Eq}}(\oI_1)-\Psi^{{\rm NEq}}(\oI^e_1).
\end{equation*}
and, hence, that the constitutive relation (\ref{Legendre-1piola-v0}) can be written more explicitly as
\begin{align}
\left\{\begin{array}{l}\bfS(\bfX,t)=2J^{-2/3}\Psi^{{\rm Eq}}_{\oI_1}\bfF+2J^{-2/3}\Psi^{{\rm NEq}}_{\oI^e_1}\bfF{\bfC^v}^{-1}-\vspace{0.2cm}\\
\hspace{1.5cm}\dfrac{2}{3}J^{-2/3}\left(I_1\Psi^{{\rm Eq}}_{\oI_1}+I^e_1\Psi^{{\rm NEq}}_{\oI^e_1}\right)\bfF^{-T}+q J\bfF^{-T}\\ \\
{\rm with}\\ \\
J=1+\dfrac{q}{\kappa}\end{array}\right.
\label{Legendre-1piola}
\end{align}

Granted the hybrid representation (\ref{Legendre-1piola}) for the stress-deformation relation of the material, the governing equations for the response of the body can be recast as the  following initial-boundary-value problem:
\begin{equation}
\left\{\begin{array}{l}{\rm Div}\left[2J^{-2/3}\Psi^{{\rm Eq}}_{\oI_1}\nabla\bfy+2J^{-2/3}\Psi^{{\rm NEq}}_{\oI^e_1}\nabla\bfy{\bfC^v}^{-1}-\right.\vspace{0.2cm}\\
\hspace{0.6cm}\dfrac{2}{3}J^{-2/3}\left(I_1\Psi^{{\rm Eq}}_{\oI_1}+I^e_1\Psi^{{\rm NEq}}_{\oI^e_1}\right)\nabla\bfy^{-T}+\vspace{0.2cm}\\
\hspace{0.6cm}\left. q J\nabla\bfy^{-T}\right]+\bfb={\bf0}, \qquad (\bfX,t)\in\Omega_0\times[0,T]\\ \\
\det\nabla\bfy-1-\dfrac{q}{\kappa}=0,\hspace{.6cm}(\bfX,t)\in\Omega_0\times[0,T]\\ \\
\bfy(\bfX,t)=\overline{\bfy}(\bfX,t), \qquad (\bfX,t)\in\partial\Omega^{\mathcal{D}}_0\times[0,T] \\ \\
\left[2J^{-2/3}\Psi^{{\rm Eq}}_{\oI_1}\nabla\bfy+2J^{-2/3}\Psi^{{\rm NEq}}_{\oI^e_1}\nabla\bfy{\bfC^v}^{-1}-\right.\vspace{0.2cm}\\
\dfrac{2}{3}J^{-2/3}\left(I_1\Psi^{{\rm Eq}}_{\oI_1}+I^e_1\Psi^{{\rm NEq}}_{\oI^e_1}\right)\nabla\bfy^{-T}+\vspace{0.2cm}\\
\left.q J\nabla\bfy^{-T}\right]\bfN=\overline{\textbf{s}}(\bfX,t), \quad (\bfX,t)\in\partial\Omega^{\mathcal{N}}_0\times[0,T] \\ \\
\bfy(\bfX,0)=\bfX, \qquad \bfX\in\Omega_0
\end{array}\right. \label{Equilibrium-PDE-Hybrid-0}
\end{equation}
coupled with the evolution equation (\ref{Evolution-ODE}), repeated here for convenience,

\begin{equation}
\left\{\begin{array}{l}
\hspace{-0.15cm}\dot{\bfC^v}(\bfX,t)=\bfG\left(\nabla\bfy^T\nabla\bfy,\bfC^v\right), \quad (\bfX,t)\in\Omega_0\times[0,T]\\ \\
\hspace{-0.15cm}\bfC^v(\bfX,0)=\bfI,\qquad\bfX\in\Omega_0\end{array}\right. \label{Evolution-ODE-Hybrid-0}
\end{equation}
for the deformation field $\bfy(\bfX,t)$, the Legendre dual field $q(\bfX,t)$, and the internal variable $\bfC^v(\bfX,t)$.

\begin{remark}
The Legendre dual field $q$ in (\ref{Equilibrium-PDE-Hybrid-0}) is nothing more that the volumetric part (\ref{Tvol}) of the Cauchy stress tensor. In other words, the Legendre variable $q$ is equal to the Cauchy hydrostatic pressure
\begin{equation*}
q=p=\dfrac{1}{3}{\rm tr}\,\bfT.
\end{equation*}
\end{remark}

\section{The time and space discretizations}\label{Sec: Discretizations}

Next, consistent with the way in which Abaqus solves initial-boundary-value problems, we discretize the governing equations (\ref{Equilibrium-PDE-Hybrid-0})-(\ref{Evolution-ODE-Hybrid-0}) in time with FD and in space with FE. We begin in Subsection \ref{Sec: Weak} by rewriting (\ref{Equilibrium-PDE-Hybrid-0}) in weak form. In Subsections \ref{Sec: Time disc} and \ref{Sec: Space disc}, we introduce the time and space discretizations one at a time. We conclude this section by  outlining a staggered method of solution for the fully discretized governing equations.

\subsection{Weak form of the governing equations}\label{Sec: Weak}

A standard calculation shows that the weak form of the initial-boundary-value problem (\ref{Equilibrium-PDE-Hybrid-0})-(\ref{Evolution-ODE-Hybrid-0}) specializes to finding $\bfy(\bfX,t)\in\mathcal{U}$ and $q(\bfX,t)\in\mathcal{V}$ such that

\begin{equation}
\left\{\begin{array}{l} \displaystyle\int_{\Omega_0}\left[2J^{-2/3}\Psi^{{\rm Eq}}_{\oI_1}\nabla\bfy+2J^{-2/3}\Psi^{{\rm NEq}}_{\oI^e_1}\nabla\bfy{\bfC^v}^{-1}-\right.\vspace{0.2cm}\\
\hspace{0.6cm}\dfrac{2}{3}J^{-2/3}\left(I_1\Psi^{{\rm Eq}}_{\oI_1}+I^e_1\Psi^{{\rm NEq}}_{\oI^e_1}\right)\nabla\bfy^{-T}+\vspace{0.2cm}\\
\hspace{0.6cm}\left. q J\nabla\bfy^{-T}\right] \cdot \nabla{\bf{w}} \,{\rm{d}}\bfX = \\
\hspace{0.6cm}\displaystyle\int_{\Omega_0} {\bfb}\cdot{\bf{w}}  \,{\rm{d}}\bfX+ \displaystyle\int_{\partial\Omega_0^{\mathcal{N}}} \overline{\bf{s}} \cdot{\bf{w}}  \,{\rm{d}}\bfX,\\
\hspace{4.0cm}\forall{\bf{w}}\in\mathcal{U}_0,~t\in\times[0,T]\\ \\
\displaystyle\int_{\Omega_0} \left(\det \nabla\bfy-1-\dfrac{q}{\kappa}\right) r \,{\rm{d}}\bfX ,\hspace{.6cm}\forall r\in\mathcal{V},~t\in\times[0,T]
\end{array}\right. \label{Weak-Form-Hybrid}
\end{equation}
with $\bfC^v(\bfX,t)$ defined by (\ref{Evolution-ODE-Hybrid-0}), where $\mathcal{U}$ and $\mathcal{V}$ are sufficiently large sets of admissible deformation $\bfy$ and pressure $q$ fields. Similarly, $\mathcal{U}_0$ stands for a sufficiently large set of test functions $\textbf{w}$. Formally, $\mathcal{U}=\{\bfy: \bfy(\bfX,t)=\overline{\bfy}(\bfX,t),\;\bfX\in\partial\Omega^{\mathcal{D}}_0\}$ and $\mathcal{U}_0=\{\bfw: \bfw(\bfX,t)=$ $\textbf{0},\;\bfX\in\partial\Omega^{\mathcal{D}}_0\}$.

\subsection{Time discretization}\label{Sec: Time disc}

Consider now a partition of the time interval $[0,T]$ into discrete times $\{t_k\}_{k=0,1,...,M}$, with $t_0=0$ and $t_M=T$. Making use of the notation
$\bfy_k(\bfX) = \bfy(\bfX,t_k)$, $q_k(\bfX) = q(\bfX,t_k)$, $\bfC^v_k(\bfX) = \bfC^v(\bfX,t_k)$, and similarly for any other time-dependent field, the governing equations (\ref{Weak-Form-Hybrid}) with (\ref{Evolution-ODE-Hybrid-0}) at any given discrete time $t_k$ take then the form

\begin{equation}
\left\{\begin{array}{l} \displaystyle\int_{\Omega_0}\left[2J_k^{-2/3}{\Psi^{{\rm Eq}}_{\oI_1}}_k\nabla\bfy_k+2J_k^{-2/3}\Psi^{{\rm NEq}}_{\oI^e_1}{}_k\nabla\bfy_k {\bfC_k^v}^{-1}-\right.\vspace{0.2cm}\\
\hspace{0.6cm}\dfrac{2}{3}J^{-2/3}_k\left({I_1}_k{\Psi^{{\rm Eq}}_{\oI_1}}_k+{I^e_1}_k{\Psi^{{\rm NEq}}_{\oI^e_1}}_k\right)\nabla\bfy_k^{-T}+\vspace{0.2cm}\\
\hspace{0.6cm}\left. q_k J_k\nabla\bfy_k^{-T}\right] \cdot \nabla{\bf{w}} \,{\rm{d}}\bfX = \\
\hspace{0.6cm}\displaystyle\int_{\Omega_0} {\bfb}_k\cdot{\bf{w}}  \,{\rm{d}}\bfX+ \displaystyle\int_{\partial\Omega_0^{\mathcal{N}}} \overline{\bf{s}}_k \cdot{\bf{w}}  \,{\rm{d}}\bfX,\quad\forall{\bf{w}}\in\mathcal{U}_0\\ \\
\displaystyle\int_{\Omega_0} \left( \det\nabla\bfy_k-1-\dfrac{q_k}{\kappa}\right) r \,{\rm{d}}\bfX=0 ,\hspace{0.8cm}\forall r\in\mathcal{V}\end{array}\right. \label{Discretized-Weak-Hybrid}
\end{equation}
and
\begin{equation}
\begin{array}{l}
\dot{\bfC^v_k}(\bfX)=\bfG(\nabla\bfy^T_k(\bfX)\nabla\bfy_k(\bfX),\bfC_k^v(\bfX)), \quad \bfX\in\Omega_0,
\end{array} \label{Discretized-ODE-Hybrid}
\end{equation}
where we emphasize that we are yet to spell out a specific time discretization for the time derivative $\dot{\bfC^v_k}$ in terms of ${\bfC^v_k}$.

\subsection{Space discretization}\label{Sec: Space disc}

Consider next a partition $^h\Omega_0 = \bigcup_{e=1}^{\texttt{N}_e} \mathcal{E}^{(e)}$ of the initial configuration $\Omega_0$ that comprises $\texttt{N}_e$ non-overlapping finite elements $\mathcal{E}^{(e)}$. Given this partition, we look for approximate solutions $^h \mathbf{y}_k(\bfX)$ and $^h q_k(\bfX)$ for the deformation $\mathbf{y}_k(\bfX)$ and pressure $q_k(\bfX)$ fields at time $t_k$ in finite dimensional subspaces of conforming finite elements. It follows that $^h\mathbf{y}_k (\bfX)$ and $^hq_k(\bfX)$ admit the representations
\begin{equation}\label{FE-Representation-y}
^h\mathbf{y}_k(\bfX) = \sum_{n=1}^{\texttt{N}_n} {}^h N^{(n)}(\bfX)\,  \mathbf{y}_k^{(n)}
\end{equation}
and
\begin{equation}\label{FE-Representation-q}
^h q_k(\bfX) = \sum_{l=1}^{{\texttt{N}_q}} {}^h N^{(l)}_{q}(\bfX) \, q_k^{(l)}
\end{equation}
in terms of the global degrees of freedom $\mathbf{y}_k^{(m)}$ and $q_k^{(l)}$, at time $t_k$, and the associated shape functions ${}^h N^{(n)}(\bfX)$ and ${}^h N^{(l)}_{q}(\bfX)$ that result from the assembly process. In these last expressions, $\texttt{N}_n$ and ${\texttt{N}_q}$ stand for the total number of nodes in the partition $^h\Omega_0$ and the total number of degrees of freedom for the approximation ${}^h q_k(\bfX)$ of the pressure field, respectively.

\begin{remark}
The shape functions ${}^h N^{(n)}(\bfX)$ and ${}^h N^{(l)}_{q}(\bfX)$ in (\ref{FE-Representation-y})-(\ref{FE-Representation-q}) must be appropriately selected so that they lead to a stable formulation and hence to a scheme that can generate converging solutions. In two recent contributions, we have made use of a class of Crouzeix-Raviart conforming finite elements \cite{GSKLP21}, as well as of a choice of elements with bubble functions \cite{WLPM23} determined from the so-called variational multiscale method \cite{Hughes95,Masud13}, that lead indeed to stable formulations. In this work, we make use of the hybrid elements built in Abaqus.
\end{remark}

By making use of the representations (\ref{FE-Representation-y})-(\ref{FE-Representation-q}), and analogous ones for the test functions $\bfw$ and $r$, equations (\ref{Discretized-Weak-Hybrid}) reduce to a system of nonlinear algebraic equations for the degrees of freedom $\mathbf{y}_k^{(m)}$ and $q_k^{(l)}$ that depend on the values, say
$^h{\bfC^v_k}$, of the internal variable $\bfC^v_k$  at the Gaussian quadrature points employed to carry out the integrals in (\ref{Discretized-Weak-Hybrid}). We write this system as
\begin{equation}
\mathcal{G}_1 \left(  {}^h\bfy_k, {}^hq_k, {}^h{\bfC^v_k} \right) = 0 .\label{Algebraic-Equilibrium}
\end{equation}
Similarly, we write the system of nonlinear algebraic equations that results from (\ref{Discretized-ODE-Hybrid}) for the internal variable ${}^h{\bfC^v_k}$ at the Gaussian quadrature points as
\begin{equation}
\mathcal{G}_2 \left(  {}^h\bfy_k , {}^h{\bfC^v_k},  {}^h\dot{\bfC}^v_k \right)= 0.\label{Algebraic-Evolution}
\end{equation}

\subsection{The solver: a staggered scheme with a fifth-order explicit Runge-Kutta integrator}\label{Sec: Solver}

Having discretized the governing equations (\ref{Equilibrium-PDE-Hybrid-0})-(\ref{Evolution-ODE-Hybrid-0}) into the system of coupled nonlinear algebraic equations (\ref{Algebraic-Equilibrium})-(\ref{Algebraic-Evolution}) for the global degrees of freedom $\bfy_k^{(m)}$, $ q_k^{(n)}$, and the internal variable $^h{\bfC^v_k}$ at the Gaussian quadrature points at time $t_k$, the final step is to solve these for given stored-energy functions $\Psi^{{\rm Eq}}$, $\Psi^{{\rm NEq}}$, given initial bulk modulus $\kappa$, given viscosity function $\eta$, given boundary data $\overline{\bfy}$ and $\overline{\bf{s}}$, and given body force $\bfb$.

Consistent with the Newton-Raphson-type solvers employed by Abaqus, we consider a staggered scheme, one in which at every time step $t_k$ the discretized equations (\ref{Algebraic-Equilibrium}) and (\ref{Algebraic-Evolution}) are solved iteratively until convergence is reached. The algorithm goes as follows:

\begin{itemize}
\item{\emph{Step 0.} Set $r=1$ and define appropriate tolerances $TOL_1$ and $TOL_2$. For a given solution ${^{h}{\bfy}_{k-1}}$, ${^{h}{q}_{k-1}}$, and ${^{h}\bfC^v_{k-1}}$ at time $t_{k-1}$, define also ${}^{h}{\bfy}_{k,r-1}={}^{h}{\bfy}_{k-1}$, ${}^{h}q_{k,r-1}={}^{h}q_{k-1}$, and ${^{h}\bfC^v_{k,r-1}}={^{h}\bfC^v_{k-1}}$.
}

\item{\emph{Step 1.} Given the boundary data $\overline{\bfy}$, $\overline{\bf{s}}$, and body force $\bfb$ at $t_k$, make use of one iteration within a Newton-Raphson solver to find ${}^{h}{\bfy}_{k,r}$ and ${}^{h}{q}_{k,r}$ such that
\begin{align}
\mathcal{G}_1 \left({}^{h}{\bfy}_{k,r}, {}^{h}{q}_{k,r},{^{h}\bfC^v_{k,r-1}}\right)&=0.
\label{SubProblem1}
\end{align}
}

\item{\emph{Step 2.} Having determined ${}^{h}{\bfy}_{k,r}$ and ${}^{h}{q}_{k,r}$ from the sub-problem (\ref{SubProblem1}), find ${^{h}\bfC^v_{k,r}}$ such that
\begin{align}
\mathcal{G}_2 \left({}^{h}{\bfy}_{k,r}, {^{h}\bfC^v_{k,r}},{^{h}\dot{\bfC}^v_{k,r}}\right)&=0.
\label{SubProblem2}
\end{align}
}

\item{\emph{Step 3.} If $\parallel \mathcal{G}_1 ({}^{h}{\bfy}_{k,r}, {}^{h}{q}_{k,r},{^{h}\bfC^v_{k,r}}) \parallel / \parallel \mathcal{G}_1 ({}^{h}{\bfy}_{k,0},$ ${}^{h}{q}_{k,0},{^{h}\bfC^v_{k,0}}) \parallel \le  TOL_1 $ and $\parallel  \mathcal{G}_2 ({}^{h}{\bfy}_{k,r},{^{h}\bfC^v_{k,r}},{^{h}\dot{\bfC}^v_{k,r}})  \parallel /  \parallel \mathcal{G}_2 ({}^{h}{\bfy}_{k,0},{^{h}\bfC^v_{k,0}},{^{h}\dot{\bfC}^v_{k,0}}) \parallel \le  TOL_2 $, then set ${^{h}\bfy_{k}}={^{h}\bfy_{k,r}}$, ${^{h}q_{k}}={^{h}q_{k,r}}$, ${^{h}\bfC^v_{k}}={^{h}\bfC^v_{k,r}}$, and move to the next time step $t_{k+1}$; otherwise set $r \leftarrow r+1$ and go back to Step 1.

}

\end{itemize}

\paragraph{The sub-problem (\ref{SubProblem2})}
The sub-problem (\ref{SubProblem2}) corresponds to a nonlinear system of first-order ODEs wherein the incompressibility constraint $\det {^{h}\bfC^v_{k,r}}=1$ is built-in. Because of the requirement of satisfying this nonlinear constraint along the entire time domain, as already remarked in the Introduction, extreme care must be exercised in the choice of the time-integration scheme. In this work, following in the footsteps of \cite{KLP16} and \cite{WLPM23}, we make use of the explicit fifth-order Runge-Kutta scheme
\begin{equation}\label{RK-5-Cv}
{^{h}\bfC^v_{k,r}}=\dfrac{1}{\left(\det\bfA_{k-1,r}\right)^{1/3}}\bfA_{k-1,r}
\end{equation}
with
\begin{multline}\label{RK-5-A}
\bfA_{k-1,r}={^{h}\bfC^v_{k-1}}+\frac{\Delta t_{k}}{90}\left(7\textbf{G}_1 + 32\textbf{G}_3 + \right.\\
\left. 12\textbf{G}_4 + 32\textbf{G}_5 + 7\textbf{G}_6\right),
\end{multline}
where
\begin{align*}
\textbf{G}_{1}=&\textbf{G}\left(\nabla{^{h}\bfy_{k-1}}, {^{h}\bfC^v_{k-1}}\right), \nonumber \\
\textbf{G}_{2}=&\textbf{G}\left(\frac{1}{2}\nabla{^{h}\bfy_{k-1}}+\frac{1}{2}\nabla{^{h}\bfy_{k,r}}, {^{h}\bfC^v_{k-1}}+ \textbf{G}_{1}\,  \frac{\mathrm{\Delta} t_k}{2}\right), \nonumber \\
\textbf{G}_{3}=&\textbf{G}\left(\frac{3}{4}\nabla{^{h}\bfy_{k-1}}+\frac{1}{4}\nabla{^{h}\bfy_{k,r}}, {^{h}\bfC^v_{k-1}}+  (3\textbf{G}_{1}+ \textbf{G}_{2})\, \frac{\mathrm{\Delta} t_k}{16}\right), \nonumber \\
\textbf{G}_{4}=&\textbf{G}\left(\frac{1}{2}\nabla{^{h}\bfy_{k-1}}+\frac{1}{2}\nabla{^{h}\bfy_{k,r}}, {^{h}\bfC^v_{k-1}}+  \textbf{G}_{3}\, \frac{\mathrm{\Delta} t_k}{2}\right), \nonumber\\
\textbf{G}_{5}=&\textbf{G}\left(\frac{1}{4}\nabla{^{h}\bfy_{k-1}}+\frac{3}{4}\nabla{^{h}\bfy_{k,r}}, {^{h}\bfC^v_{k-1}}+ \right.\\
&\left. 3 (-\textbf{G}_{2}+ 2\textbf{G}_{3}+ 3\textbf{G}_{4})\, \frac{\mathrm{\Delta} t_k}{16}\right),  \\
\textbf{G}_{6}=&\textbf{G}\left(\nabla{^{h}\bfy_{k,r}}, {^{h}\bfC^v_{k-1}}+  (\textbf{G}_{1}+ 4\textbf{G}_{2} +\right.\\
&\left.  6\textbf{G}_{3}- 12\textbf{G}_{4}+ 8\textbf{G}_{5})\, \frac{\mathrm{\Delta} t_k}{7}\right),  %\label{Runge-Kutta-k}
\end{align*}
$\mathrm{\Delta} t_k=t_{k}-t_{k-1}$, and where we recall that the function $\textbf{G}$ is defined by (\ref{Evolution-I1-J}).

\begin{remark}
The scheme (\ref{RK-5-Cv})-(\ref{RK-5-A}) is a modification of a scheme originally introduced by Lawson \cite{Lawson66} that has the distinctive merit of preserving the constraint $\det {^{h}\bfC^v_{k,r}}=1$ identically \cite{WLPM23}. Recent numerical experiments \cite{GSKLP21,SGLP23,WLPM23} have shown that (\ref{RK-5-Cv})-(\ref{RK-5-A}) is a highly performant scheme, thus its use here.
\end{remark}
\begin{remark}
The time increment $\mathrm{\Delta} t_k=t_{k}-t_{k-1}$ in (\ref{RK-5-Cv})-(\ref{RK-5-A}) must be chosen to be sufficiently small relative to the material time scale (\ref{tau}) in the evolution equation (\ref{Evolution-I1-J}). In practice, it suffices to set $\Delta t_k \leq 10^{-2}\tau$. Note that for the explicit scheme used here, the choice of $\Delta t_k$ indirectly sets the value of $TOL_2$ in \emph{Step 3} above. For other schemes, such as implicit and/or adaptive schemes, the value of $\Delta t_k$ could be adjusted in terms of a desired $TOL_2$.
\end{remark}

\section{The UMAT implementation}\label{Sec: UMAT}

The discretized governing equations (\ref{SubProblem1}) and (\ref{SubProblem2}) can be assembled and solved in Abaqus by making use of a UMAT subroutine in conjunction with hybrid finite elements.

In particular, on one hand, the Runge-Kutta scheme (\ref{RK-5-Cv})-(\ref{RK-5-A}) can be directly coded within the UMAT subroutine to solve the equations (\ref{SubProblem2}) for the values ${^{h}\bfC^v_{k,r}}$ of the internal variable  ${\bfC^v_{k}}$ at the Gaussian quadrature points at any time increment $t_k$ and Newton-Raphson iteration $r$.

On the other hand, to solve the equations (\ref{SubProblem1}) for the degrees of freedom $\bfy_k^{(n)}$ and $q_k^{(l)}$ describing the FE approximations ${}^{h}\bfy_k(\bfX)$ and ${}^{h}q_k(\bfX)$ for the deformation and pressure fields at the time increment $t_k$, Abaqus requires to code within the UMAT subroutine the stress-deformation relation of the material and its derivatives with respect to the deformation in a certain format. Specifically, one must define the Cauchy stress
\begin{multline*}
\bfT=2J^{-5/3}\Psi^{{\rm Eq}}_{\oI_1}\bfF\bfF^T+2J^{-5/3}\Psi^{{\rm NEq}}_{\oI^e_1}\bfF{\bfC^v}^{-1}\bfF^T-\\
\left(\dfrac{2}{3}J^{-5/3}I_1\Psi^{{\rm Eq}}_{\oI_1}+\dfrac{2}{3}J^{-5/3}I^e_1\Psi^{{\rm NEq}}_{\oI^e_1}\right)\bfI+\kappa(\widehat{J}-1)\bfI
\end{multline*}
in terms of the deformation gradient tensor $\bfF$, the variable

\begin{equation*}
\widehat{J}=1+\dfrac{q}{\kappa},
\end{equation*}
and the internal variable $\bfC^v$. For the derivatives of the stress-deformation relation, one must define the volumetric moduli
\begin{equation*}
\widehat{K}=J\dfrac{\partial^2}{\partial{\widehat{J}}^2}\left[\dfrac{\kappa}{2}(\widehat{J}-1)^2\right]=\kappa J \quad {\rm and}\quad \dfrac{\partial \widehat{K}}{\partial \widehat{J}}=0,
\end{equation*}
as well as the tangent modulus
\begin{align*}
\mathcal{L}_{ijkl}=&\dfrac{1}{2J}\left(\dfrac{\partial ({\rm dev}\,\tau_{ij})}{\partial F_{kr}}F_{lr}+\dfrac{\partial  ({\rm dev}\,\tau_{ij})}{\partial F_{lr}}F_{kr}\right)+\kappa J\delta_{ij}\delta_{kl}\\
=&4J^{-1}\Psi^{{\rm Eq}}_{\oI_1\oI_1}\left(\overline{B}_{ij}-\dfrac{\oI_1}{3}\delta_{ij}\right)\left(\overline{B}_{kl}-\dfrac{\oI_1}{3}\delta_{kl}\right)+\\
&4J^{-1}\Psi^{{\rm NEq}}_{\oI^e_1\oI^e_1}\left(\overline{B}^e_{ij}-\dfrac{\oI^e_1}{3}\delta_{ij}\right)\left(\overline{B}^e_{kl}-\dfrac{\oI^e_1}{3}\delta_{kl}\right)+\\
&2J^{-1}\Psi^{{\rm Eq}}_{\oI_1}\left[\dfrac{1}{2}\left(\delta_{ik}\overline{B}_{jl}+\delta_{jk}\overline{B}_{il}+\delta_{il}\overline{B}_{jk}+\delta_{jl}\overline{B}_{ik}\right)-\right.\\
&\left.\dfrac{2}{3}\left(\overline{B}_{ij}\delta_{kl}+\delta_{ij}\overline{B}_{kl}\right)+\dfrac{2}{9}\oI_1\delta_{ij}\delta_{kl}\right]\\
&2J^{-1}\Psi^{{\rm NEq}}_{\oI^e_1}\left[\dfrac{1}{2}\left(\delta_{ik}\overline{B}^e_{jl}+\delta_{jk}\overline{B}^e_{il}+\delta_{il}\overline{B}^e_{jk}+\delta_{jl}\overline{B}^e_{ik}\right)-\right.\\
&\left.\dfrac{2}{3}\left(\overline{B}^e_{ij}\delta_{kl}+\delta_{ij}\overline{B}^e_{kl}\right)+\dfrac{2}{9}\oI^e_1\delta_{ij}\delta_{kl}\right]+\kappa J\delta_{ij}\delta_{kl}
\end{align*}
in terms of the deviatoric part of the Kirchhoff stress tensor $\boldsymbol{\tau}=J\bfT$ and the alternate bulk modulus $\widehat{K}=\kappa J$, where, for simplicity, we have made use of the notation $\obfB=J^{-2/3}\bfF\bfF^T$ and $\obfB^e=J^{-2/3}\bfF{\bfC^v}^{-1}\bfF^T$.

As an example that contains all nonlinearities, we have implemented the UMAT subroutine for the case when the equilibrium $\Psi^{{\rm Eq}}$ and non-equilibrium $\Psi^{{\rm NEq}}$ stored-energy functions are given by (\ref{WLP-Eq}) and (\ref{WLP-NEq}) and the viscosity $\eta$ is given by (\ref{eta-KLP}). The subroutine is available in GitHub.\footnote{\url{https://github.com/victorlefevre/UMAT_Lefevre_Sozio_Lopez-Pamies}}

\section{Error analysis}\label{Sec: Convergence}

In this section, we demonstrate the accuracy and convergence properties of the proposed UMAT subroutine by first performing a patch test, where the fields are homogeneous, and then comparing with an exact solution wherein the fields are non-homogeneous. Beyond demonstrating its accuracy and convergence, the results serve to illustrate that the subroutine can be used with any of the hybrid simplicial and hexahedral elements built in Abaqus for 3D problems, as well as with its hybrid 2D axisymmetric quadrilateral elements.

All the results that are presented in this section pertain to the case when the equilibrium $\Psi^{{\rm Eq}}$ and non-equilibrium $\Psi^{{\rm NEq}}$ stored-energy functions are given by (\ref{WLP-Eq}) and (\ref{WLP-NEq}) and the viscosity $\eta$ is given by (\ref{eta-KLP}) with the material constants listed in Table \ref{Table1}. This choice of constitutive model and materials constants is descriptive of the acrylic elastomer VHB 4910 from 3M \cite{Steinmann12,KLP16}, which is strongly nonlinear in elasticity and  viscosity and hence an ideal test case to probe the accuracy and convergence of the subroutine.
\begin{table}[h!]\centering
\caption{Values of the material constants in the stored-energy and viscosity functions (\ref{WLP-Eq})-(\ref{eta-KLP}) for the acrylic elastomer VHB 4910.}
\begin{tabular}{lll}
\toprule
$\mu_1=\SI{13.54}{\kilo\pascal}$  & $m_1=\SI{5.42}{\kilo\pascal}$ & $\eta_0=\SI{7014}{\kilo\pascal\second}$ \\
$\alpha_1=1.00$  & $a_1=-10$ & $\eta_{\infty}=\SI{0.1}{\kilo\pascal\second}$\\
$\mu_2=\SI{1.08}{\kilo\pascal}$  & $m_2=\SI{20.78}{\kilo\pascal}$  &  $K_1=\SI{3507}{\kilo\pascal\second}$    \\
$\alpha_2=-2.474$ & $a_2=1.948$ & $K_2=\SI{1}{\kilo\pascal^{-2}}$   \\
& &    $\beta_1=1.852$     \\
& &    $\beta_2=0.26$      \\
\bottomrule
\end{tabular} \label{Table1}
\end{table}

%
%\begin{table}[h!]\centering
%\caption{Values of the material constants in the stored-energy and viscosity functions (\ref{WLP-Eq})-(\ref{eta-KLP}) for the acrylic elastomer VHB 4910.}
%\begin{tabular}{llll}
%\toprule
%$\mu_1=13.54$ kPa & & & $\mu_2=1.08$ kPa  \\
%$\alpha_1=1.00$ & & & $\alpha_2=-2.474$  \\
%$m_1=5.42$ kPa & & & $m_2=20.78$ kPa \\
%$a_1=-10$ & & & $a_2=1.948$ \\
%$\eta_0=7014$ kPa s & & & $\eta_{\infty}=0.1$ kPa s \\
%$K_1=3507$ kPa s & & & $K_2=1$ kPa$^{-2}$  \\
%$\beta_1=1.852$ & & & $\beta_2=0.26$ \\
%\bottomrule
%\end{tabular} \label{Table1}
%\end{table}
%

\subsection{Patch test}

For the patch test, we consider a unit cube occupying the domain
\begin{equation*}
\Omega_0=\{\bfX: 0<X_1<1,\,0<X_2<1,\,0<X_3<1\}
\end{equation*}
with respect to the Cartesian laboratory frame of reference $\{\bfe_1,\bfe_2,\bfe_3\}$ that is subjected to loading/unloading in uniaxial tension at a constant stretch rate in the $\bfe_3$ direction. Precisely, we set the body force to $\textbf{b}=\textbf{0}$ and consider that the cube is subjected to the affine boundary conditions
\begin{equation*}
\left\{\begin{array}{ll}
y_3(\bfX,t)=0, & \bfX=X_1\bfe_1+X_2\bfe_2\vspace{0.2cm}\\
y_3(\bfX,t)=\overline{F}_{33}(t), &\bfX=X_1\bfe_1+X_2\bfe_2+\bfe_3
\end{array}\right. ,   %\label{BCs-cube}
\end{equation*}
%
%
%%%%%%%%%%%%%%%%%%%%%%%%%%%%%%%%%%%%%%%%%%%%%%%%%%%%%%%%%%%%%%%%%%%%%%%%%%%%%%
\begin{figure}[t!]
  \subfigure[]{
   \label{fig:2a}
   \begin{minipage}[]{0.5\textwidth}
   \centering \includegraphics[width=2.55in]{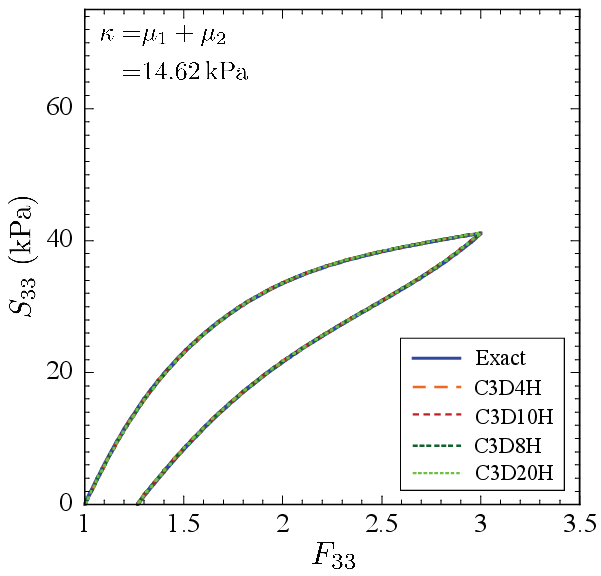}
   \vspace{0.2cm}
   \end{minipage}}
  \subfigure[]{
   \label{fig:2b}
   \begin{minipage}[]{0.5\textwidth}
   \centering \includegraphics[width=2.55in]{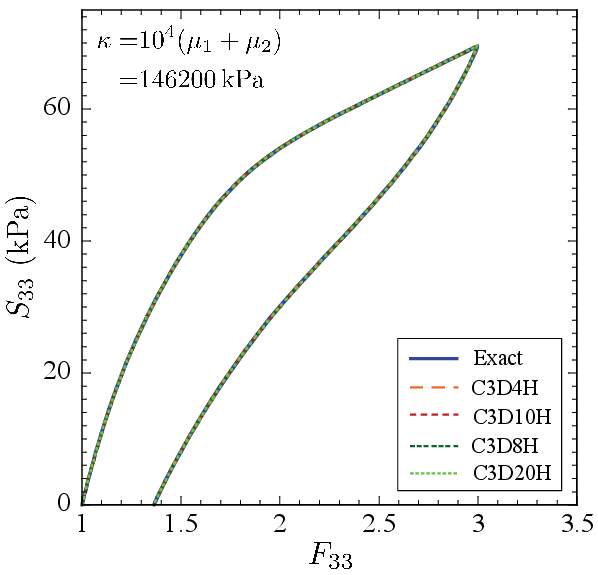}
   \vspace{0.2cm}
   \end{minipage}}
   \caption{Comparison between the FE results obtained with the Abaqus hybrid elements C3D4H, C3D10H, C3D8H, C3D20H and the exact solution for the response of (a) a highly compressible and (b) a nearly incompressible viscoelastic elastomer --- with stored-energy functions (\ref{WLP-Eq}), (\ref{WLP-NEq}), viscosity (\ref{eta-KLP}), and the material constants listed in Table \ref{Table1} --- under uniaxial tension loading/unloading at the constant stretch rate $|\dot{\lambda}_0|=0.05$ s$^{-1}$.}\label{Fig2}
\end{figure}
%%%%%%%%%%%%%%%%%%%%%%%%%%%%%%%%%%%%%%%%%%%%%%%%%%%%%%%%%%%%%%%%%%%%%%%%%%%%%%
%
%
\begin{equation*}
\left\{\begin{array}{ll}
\bfe_1\cdot\bfS(\bfX,t)\bfe_3=0, & \bfX=X_1\bfe_1+X_2\bfe_2\vspace{0.2cm}\\
\bfe_2\cdot\bfS(\bfX,t)\bfe_3=0, & \bfX=X_1\bfe_1+X_2\bfe_2\vspace{0.2cm}\\
\bfe_1\cdot\bfS(\bfX,t)\bfe_3=0, & \bfX=X_1\bfe_1+X_2\bfe_2+\bfe_3\vspace{0.2cm}\\
\bfe_2\cdot\bfS(\bfX,t)\bfe_3=0, & \bfX=X_1\bfe_1+X_2\bfe_2+\bfe_3\vspace{0.2cm}\\
\end{array}\right. ,   %\label{BCs-cube}
\end{equation*}
and
\begin{equation*}
\left\{\begin{array}{ll}
\bfS(\bfX,t)\bfe_1=\textbf{0}, & \bfX=X_2\bfe_2+X_3\bfe_3\vspace{0.2cm}\\
\bfS(\bfX,t)\bfe_1=\textbf{0}, & \bfX=\bfe_1+X_2\bfe_2+X_3\bfe_3\vspace{0.2cm}\\
\bfS(\bfX,t)\bfe_2=\textbf{0}, &\bfX=X_1\bfe_1+X_3\bfe_3\vspace{0.2cm}\\
\bfS(\bfX,t)\bfe_2=\textbf{0}, & \bfX=X_1\bfe_1+\bfe_2+X_3\bfe_3
\end{array}\right. ,  %\label{BCs-cube}
\end{equation*}
where
\begin{equation*}
\overline{F}_{33}(t)=\left\{\begin{array}{ll}
1+\dot{\lambda}_0 t, & 0\leq t\leq \dfrac{2}{\dot{\lambda}_0}\vspace{0.2cm}\\
5-\dot{\lambda}_0 t, & \dfrac{2}{\dot{\lambda}_0}<t\leq \dfrac{4}{\dot{\lambda}_0}\vspace{0.2cm}\end{array}\right. {\rm with}\; \dot{\lambda}_0=0.05\,{\rm s}^{-1}.
\end{equation*}

Figure \ref{Fig2} compares with the exact solution the results for the stress component $S_{33}$ in the cube as computed with the Abaqus hybrid simplicial elements C3D4H and C3D10H and the hexahedral elements C3D8H and C3D20H. The results in Fig. \ref{Fig2}(a) correspond to a highly compressible elastomer with initial bulk modulus $\kappa=\mu_1+\mu_2=14.62$ kPa, while those in Fig. \ref{Fig2}(b) correspond to a nearly incompressible elastomer with initial bulk modulus $\kappa=10^4(\mu_1+\mu_2)=146200$ kPa.

It is plain from the comparisons that the code passes the patch test for all four types of elements and both compressibilities. Many other different types of patch tests (not included here for conciseness of presentation) have corroborated the satisfaction of this basic completeness requirement for the code.

\subsection{The radially symmetric deformation of a spherical shell}

In the sequel, we analyze the accuracy and convergence of the solutions generated by the subroutine as the size $h$ of the finite elements decreases --- that is, as the total number $\texttt{N}_{dof}$ of degrees of freedom increases ---  by direct comparison with the exact solution for one of the few initial-boundary-value problems in finite viscoelastostatics for which an exact solution can be determined in closed form, that of the radially symmetric deformation of a spherical shell made of an incompressible material \cite{KAILP17}.

In particular, we consider a spherical shell of initial inner radius $A=0.9$ m and initial outer radius $B=1$ m. At time $t=0$, with respect to the Cartesian laboratory
frame of reference $\{\bfe_1,\bfe_2,\bfe_3\}$, the shell occupies the domain
\begin{equation*}
\Omega_0=\{\bfX: A<R<B\}\quad{\rm with}\quad R\equiv|\bfX|=\sqrt{\bfX\cdot\bfX}.
\end{equation*}

Again, much like for the preceding patch test, we take the shell to be made of an elastomer with non-Gaussian stored-energy functions (\ref{WLP-Eq}) and (\ref{WLP-NEq}), and nonlinear viscosity (\ref{eta-KLP}), with the material constants listed in Table \ref{Table1} for the acrylic elastomer VHB 4910. We restrict attention to the limiting case when the elastomer is fully incompressible and hence set $\kappa=+\infty$.

The inner boundary of the shell is traction free, while its outer boundary is subjected to a prescribed radial deformation. Precisely, the boundary conditions are given by
\begin{equation}
\left\{\begin{array}{ll}
\bfy(\bfX,t)=\dfrac{b(t)}{B}\bfX, & (\bfX,t)\in\partial\Omega_0^{\mathcal{D}}\times [0,T]\\ \\
-\bfS(\bfX,t)\bfX=0, & (\bfX,t)\in\partial\Omega_0^{\mathcal{N}}\times [0,T]\end{array}\right. , \label{BCs-shell}
\end{equation}
where $b(t)$ is the prescribed value of the outer radius $r=|\bfx|$ of the shell at time $t$ and
\begin{equation*}
\partial\Omega_0^{\mathcal{D}}=\{\bfX: R=B\} \quad {\rm and}\quad \partial\Omega_0^{\mathcal{N}}=\{\bfX: R=A\}.
\end{equation*}

\subsubsection{The exact solution}

Restricting attention to radially symmetric solutions, it follows from the incompressibility constraint $\det\nabla\bfy=1$ and the Dirichlet boundary condition (\ref{BCs-shell})$_1$ that the deformation field is given by the fully explicit relation \cite{KAILP17}
\begin{equation*}
\bfy(\bfX,t)=\lambda(R,t)\bfX
\end{equation*}
with
\begin{equation*}%\label{l-Shell}
\lambda(R,t)=\left(1+\dfrac{b(t)^3-B^3}{R^3}\right)^{1/3}.
\end{equation*}
In turn, absent body forces, it can be shown from the remaining governing equations that the first Piola-Kirchhoff stress is given by \cite{KAILP17}
\begin{equation}\label{S-Shell}
\bfS(\bfX,t)=s_1(R,t)\dfrac{1}{R^2}\bfX\otimes\bfX+s_2(R,t)\left(\bfI-\dfrac{1}{R^2}\bfX\otimes\bfX\right)
\end{equation}
with
\begin{align}\label{s1-Shell}
s_1(R,t)=&\lambda^2(R,t)\displaystyle\int_{A}^{R}\dfrac{1}{z\lambda^2(z,t)}\dfrac{{\rm d}\mathcal{W}^{{\rm Eq}}}{{\rm d}\lambda}(\lambda(z,t))\,{\rm d}z+\nonumber\\
&\lambda^2(R,t)\displaystyle\int_{A}^{R}\dfrac{1}{z\lambda^2(z,t)}\dfrac{\partial\mathcal{W}^{{\rm NEq}}}{\partial\lambda}(\lambda(z,t),\lambda_v(z,t))\,{\rm d}z
\end{align}
and
\begin{align*}
s_2(R,t)=&\lambda^{-3}(R,t)s_1(R,t)+\\
&2\Psi^{{\rm Eq}}_{\oI_1}(\oI_1)\left(\lambda(R,t)-\lambda^{-5}(R,t)\right)+\\
&2\Psi^{{\rm NEq}}_{\oI^e_1}(\oI^e_1)\left(\dfrac{\lambda(R,t)}{\lambda_v^{2}(R,t)}-\dfrac{\lambda_v^4(R,t)}{\lambda^{5}(R,t)}\right),
\end{align*}
where
\begin{align*}
\mathcal{W}^{{\rm Eq}}(\lambda)\equiv\Psi^{{\rm Eq}}(\oI_1),\quad \mathcal{W}^{{\rm NEq}}(\lambda,\lambda_v)\equiv\Psi^{{\rm NEq}}(\oI^e_1),
\end{align*}
$\oI_1=(1+2\lambda^6)/\lambda^4$, $\oI_1^e=(2\lambda^6+\lambda_v^6)/(\lambda^4\lambda_v^2)$, $\oI_2^e=(\lambda^6+2\lambda_v^6)/(\lambda^2\lambda_v^4)$, and
$\oI_1^v=(1+2\lambda_v^6)/\lambda_v^4$, and where $\lambda_v(R,t)$ is defined implicitly as the solution of the evolution equation
\begin{equation}
\left\{\begin{array}{l}\dot{\lambda}_v(R,t)=\dfrac{\Psi_{\oI_1^e}^{{\rm NEq}}(\oI_1^e){\lambda}_v(R,t)\left(\dfrac{\lambda^6(R,t)}{{\lambda}^6_v(R,t)}-1\right)}{3\eta(\oI_1^e,\oI_2^e,I_1^v)\dfrac{\lambda^4(R,t)}{{\lambda}^4_v(R,t)}} \vspace{0.2cm} \\
{\lambda}_v(R,0)=1\end{array}\right. .  \label{Evolution-lv}
\end{equation}

\subsubsection{The FE results versus the exact solution}

%
%%%%%%%%%%%%%%%%%%%%%%%%%%%%%%%%%%%%%%%%%%%%%%%%%%%%%%%%%%%%%%%%%%%%%%%%%%%%%%
\begin{figure}[t!]
  \subfigure[]{
   \label{fig:2a}
   \begin{minipage}[]{0.5\textwidth}
   \centering \includegraphics[width=2.5in]{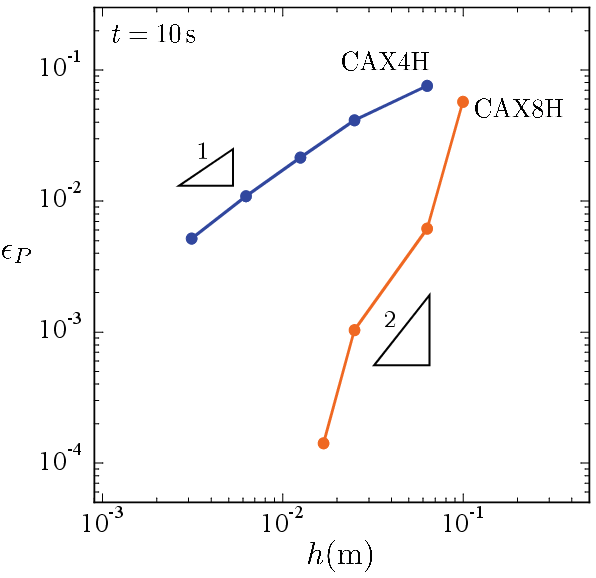}
   \vspace{0.2cm}
   \end{minipage}}
  \subfigure[]{
   \label{fig:2b}
   \begin{minipage}[]{0.5\textwidth}
   \centering \includegraphics[width=2.5in]{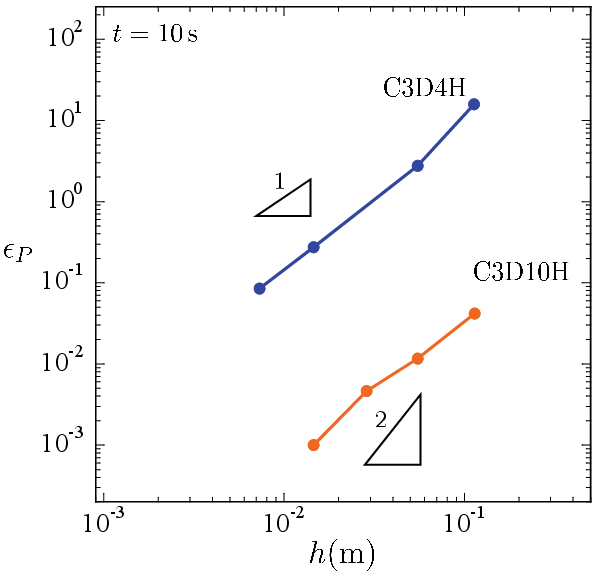}
   \vspace{0.2cm}
   \end{minipage}}
     \subfigure[]{
   \label{fig:2b}
   \begin{minipage}[]{0.5\textwidth}
   \centering \includegraphics[width=2.5in]{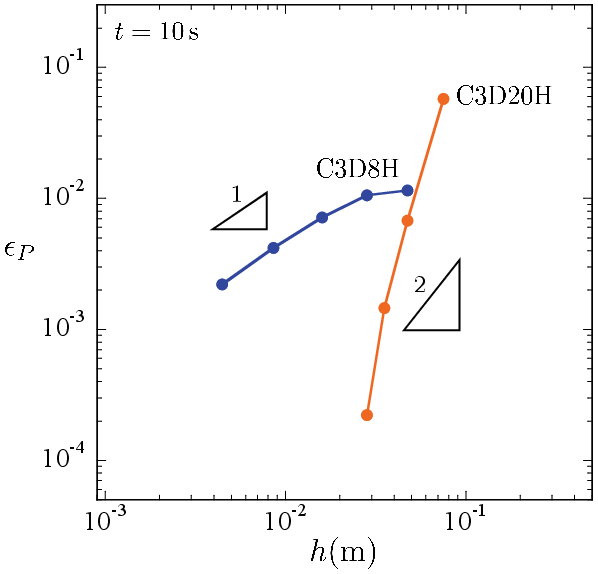}
   \vspace{0.2cm}
   \end{minipage}}
   \caption{The error (\ref{Error-P}) at time $t=10$ s as a function of the average mesh size $h$ for: (a) 2D quadrilateral axisymmetric elements CAX4H and CAX8H, (b) 3D simplicial elements C3D4H and C3D10H, and (c) hexahedral elements C3D8H and C3D20H.}\label{Fig3}
\end{figure}
%%%%%%%%%%%%%%%%%%%%%%%%%%%%%%%%%%%%%%%%%%%%%%%%%%%%%%%%%%%%%%%%%%%%%%%%%%%%%%
%

It follows immediately from the stress relations (\ref{S-Shell}) and (\ref{s1-Shell}) that the nominal pressure at the outer boundary of the shell that results from applying the radial deformation (\ref{BCs-shell})$_1$ is given by
\begin{align}\label{P-Shell}
P(t)=&s_1(B,t)\nonumber\\
=&\dfrac{b^2(t)}{B^2}\displaystyle\int_{A}^{B}\dfrac{1}{R\lambda^2(R,t)}\dfrac{{\rm d}\mathcal{W}^{{\rm Eq}}}{{\rm d}\lambda}(\lambda(R,t))\,{\rm d}R+\nonumber\\
&\dfrac{b^2(t)}{B^2}\displaystyle\int_{A}^{B}\dfrac{1}{R\lambda^2(R,t)}\dfrac{\partial\mathcal{W}^{{\rm NEq}}}{\partial\lambda}(\lambda(R,t),\lambda_v(R,t))\,{\rm d}R.
\end{align}

In Fig. \ref{Fig3}, we confront with the exact solution (\ref{P-Shell}) the corresponding pressure $P^h(t)$ generated by the proposed UMAT subroutine for several of the Abaqus hybrid elements and a range of finite element sizes $h$ for the case when the shell is deformed at the constant deformation rate
\begin{align*}%\label{bt-Shell}
b(t)=B+B\dot{\lambda}_0 t\quad {\rm with}\quad \dot{\lambda}_0=0.05\,{\rm s}^{-1}.
\end{align*}
\begin{remark}
In order to evaluate the exact solution (\ref{P-Shell}) for the pressure $P(t)$ accurately, several approaches are possible. In this work, we employ Gaussian quadrature for the second integral in (\ref{P-Shell}), which requires solving the evolution equation (\ref{Evolution-lv}) for $\lambda_v(R_i,t)$ at each of the Gauss points $R_i$. We do so by making use of the explicit fifth-order-accurate Runge-Kutta scheme of Lawson \cite{Lawson66}. For the problem at hand here, $\texttt{N}_{Gauss}=100$ Gauss points suffice to deliver accurate values for $P(t)$.
\end{remark}

Specifically, Fig. \ref{Fig3}(a) presents results for the error
\begin{align}\label{Error-P}
\epsilon_{P}\equiv\dfrac{|P(t)-P^h(t)|}{|P(t)|}
\end{align}
between the FE approximation $P^h(t)$ and the exact solution (\ref{P-Shell}) as a function of the average mesh size $h$ for the cases when the FE result is computed with the Abaqus hybrid 2D quadrilateral axisymmetric elements CAX4H and CAX8H. The results pertain to the fixed time $t=10$ s, which corresponds to the instance at which the deformed outer radius is $b(10)=1.5$ m.

Figures \ref{Fig3}(b) and \ref{Fig3}(c) present results entirely analogous to those presented in Fig. \ref{Fig3}(a) for the cases when the FE result is computed with the Abaqus hybrid simplicial elements C3D4H and C3D10H and the hexahedral elements C3D8H and C3D20H, respectively.

It is immediate from Fig. \ref{Fig3} that the proposed UMAT subroutine generates solutions that converge to the exact solution at roughly the expected rates for all six types of elements. It is also clear from the figures that, with the exception of the low-order simplicial element C3D4H, all other types of elements lead to accurate solutions even for meshes that are only moderately refined.

\section{Application to bicontinuous rubber blends}\label{Sec: Performance}

\begin{figure*}[t!]
\centering
\includegraphics[width=0.85\textwidth]{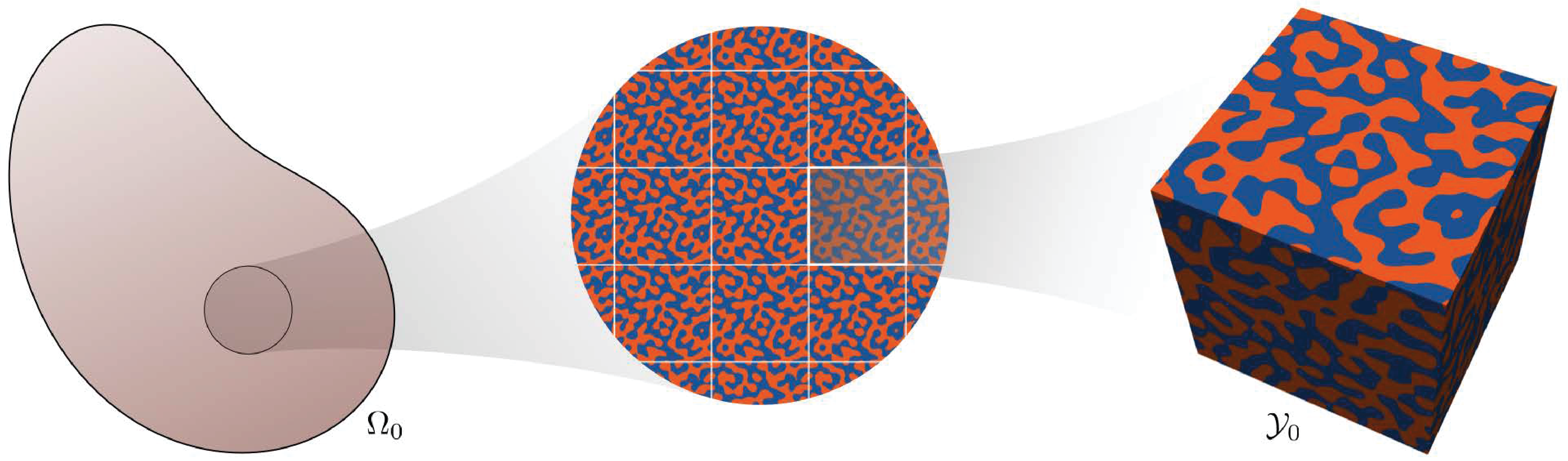}
\caption{Schematics (in the initial configuration) of a bicontinuous rubber blend and of the unit cell $\mathcal{Y}_0$ that defines its periodic microstructure.}
\label{Fig4}
\end{figure*}

Combined with the ample and well-established capabilities of Abaqus, the UMAT subroutine introduced in this work provides a powerful tool to study a vast spectrum of problems, such as, for instance, the indentation \cite{Shield52,Hu10,Cai17}, the homogenization \cite{GSKLP21,SGLP23}, and the fracture \cite{SLP23a,SLP23b,SLP23c} of viscoelastic elastomers, for which one can leverage the capabilities built in Abaqus to deal with contact, periodic boundary conditions, and cracks.

In this last section, by way of an example to showcase its capabilities, we present an application of the proposed UMAT subroutine to solve a homogenization problem, that of a bicontinuous rubber blend.

Because of their unique mechanical and physical properties, rubber blends have long been a staple in countless technological applications. A majority of the rubber blends that are utilized in applications exhibit bicontinuous microstructures, that is, they are binary mixtures in which each constituent or phase is segregated into an interpenetrating network of two separate but fully continuous domains that are perfectly bonded to one another \cite{ryan2002designer,pernot2002design}. The characteristic length scales of these microstructures is typically in the order of at most a micron and hence they are small enough that the macroscopic behavior (at the length scale of millimeters and larger) of such blends is expected to be accurately described by a homogenization limit.

\subsection{The homogenized viscoelastic response of a bicontinuous blend of Gaussian rubbers with constant viscosity}

Following in the footsteps of \cite{SLPLP24}, as illustrated schematically by Fig. \ref{Fig4}, we consider the idealization of a bicontinuous rubber blend as the periodic repetition of a unit cell $\mathcal{Y}_0$ made of two rubber phases --- labeled $r=1$ and $r=2$ --- whose initial spatial distributions at time $t=0$ are described by the characteristic or indicator functions
\begin{equation}\label{theta-blend}
\theta_0^{(r)}(\bfX)=\left\{\begin{array}{ll}1 & \textrm{if }\bfX \textrm{ is in phase } r\\
0 &\textrm{else}\end{array}\right.\quad r=1,2.
\end{equation}

For clarity of presentation, we take the two rubber phases in the blend to be isotropic and nearly incompressible viscoelastic solids with Gaussian elasticity and constant viscosity. Accordingly, the first Piola-Kirchhoff stress tensor $\bfS$ at any material point $\bfX$ and time $t$ is given by \cite{KLP16,GSKLP21}
\begin{align*}
\bfS(\bfX,t)=&\frac{\partial \psi}{\partial\bfF}(\bfX,\bfF,\bfC^v)\\
=&J^{-2/3}\mu(\bfX)\bfF+J^{-2/3} m(\bfX)\bfF{\bfC^{v}}^{-1}-\\
&\frac{J^{-2/3}}{3}\left(I_1\mu(\bfX)+I_1^e m(\bfX)\right)+\kappa(\bfX)(J-1)J\bfF^{-T},%\label{S-blend}
\end{align*}
where $\bfC^v$ is implicitly defined by the evolution equation

\begin{equation*}
\left\{\begin{array}{l}\dot{\bfC}^v(\bfX,t)=\dfrac{J^{-2/3}m(\bfX)}{\eta(\bfX)}\left[\bfC-\dfrac{1}{3}\left(\bfC\cdot{\bfC^v}^{-1}\right)\bfC^v\right] \\ \\
\bfC^v(\bfX,0)=\bfI\end{array}\right.,  %\label{Evolution-blend}
\end{equation*}
and where $\mu(\bfX)=\theta_0^{(1)}(\bfX)\mu^{(1)}+\theta_0^{(2)}(\bfX)\mu^{(2)}$, $m(\bfX)=\theta_0^{(1)}(\bfX)$ $m^{(1)}+\theta_0^{(2)}(\bfX)m^{(2)}$,  $\kappa(\bfX)=\theta_0^{(1)}(\bfX)\kappa^{(1)}+\theta_0^{(2)}(\bfX)$ $\kappa^{(2)}$, and $\eta(\bfX)=\theta_0^{(1)}(\bfX)\eta^{(1)}+\theta_0^{(2)}(\bfX)\eta^{(2)}$ stand for the pointwise equilibrium and non-equilibrium initial shear moduli, the initial bulk modulus, and the viscosity of the blend.

In this setting, as recently established in \cite{GSKLP21}, the homogenized viscoelastic response of the blend is defined by the relation between the history of the macroscopic first Piola-Kirchhoff stress tensor
\begin{equation}
\left\{\boldsymbol{\mathsf{S}}(t),\;t\in[0,T]\right\},\qquad \boldsymbol{\mathsf{S}}(t)\equiv\dfrac{1}{|\mathcal{Y}_0|}\int_{\mathcal{Y}_{0}}\bfS(\bfX,t)\,{\rm d}\bfX, \label{S-avg-Y}
\end{equation}
and the history of the macroscopic deformation gradient tensor
\begin{equation}
\left\{\boldsymbol{\mathsf{F}}(t),\;t\in[0,T]\right\},\qquad \boldsymbol{\mathsf{F}}(t)\equiv\dfrac{1}{|\mathcal{Y}_0|}\int_{\mathcal{Y}_{0}}\bfF(\bfX,t)\,{\rm d}\bfX, \label{F-avg-Y}
\end{equation}
based on the solution for the \emph{unit-cell problem}
\begin{equation}
\left\{\begin{array}{l}
\hspace{-0.15cm}{\rm Div}\left[J^{-2/3}\mu(\bfX)\nabla\bfy+J^{-2/3}m(\bfX)\nabla\bfy{\bfC^v}^{-1}-\right.\vspace{0.2cm}\\
\hspace{-0.15cm}\hspace{0.6cm}\dfrac{1}{3}J^{-2/3}\left(I_1\mu(\bfX)+I^e_1 m(\bfX)\right)\nabla\bfy^{-T}+\vspace{0.2cm}\\
\hspace{-0.15cm}\hspace{0.6cm}\left. q J\nabla\bfy^{-T}\right]={\bf0}, \qquad (\bfX,t)\in\mathcal{Y}_0\times[0,T]\\ \\
\hspace{-0.15cm}\det\nabla\bfy-1-\dfrac{q}{\kappa(\bfX)}=0,\hspace{.6cm}(\bfX,t)\in\mathcal{Y}_0\times[0,T]\\ \\
\hspace{-0.15cm}\bfy(\bfX,0)=\bfX, \qquad \bfX\in\mathcal{Y}_0
\end{array}\right. \label{Equilibrium-PDE-Hybrid-blend}
\end{equation}
coupled with the evolution equation
\begin{equation}
\left\{\begin{array}{l}\hspace{-0.15cm}\dot{\bfC}^v(\bfX,t)=\dfrac{J^{-2/3}m(\bfX)}{\eta(\bfX)}\left[\nabla\bfy^T\nabla\bfy-\right.\vspace{0.2cm}\\
\hspace{3.35cm}\left.\dfrac{1}{3}\left(\nabla\bfy^T\nabla\bfy\cdot{\bfC^v}^{-1}\right)\bfC^v\right] \\ \\
\hspace{-0.15cm}\bfC^v(\bfX,0)=\bfI\end{array}\right. \label{Evolution-ODE-Hybrid-blend}
\end{equation}
for the deformation field
\begin{align*}
\bfy(\bfX,t)=\boldsymbol{\mathsf{F}}(t)\bfX+\widetilde{\bfy}(\bfX,t),
\end{align*}
where $\boldsymbol{\mathsf{F}}(t)$ is prescribed and $\widetilde{\bfy}(\bfX,t)$ is $\mathcal{Y}_0$-periodic, the $\mathcal{Y}_0$-periodic pressure field $q(\bfX,t)$, and the internal variable $\bfC^v(\bfX,t)$.

\begin{remark}
The unit-cell problem (\ref{Equilibrium-PDE-Hybrid-blend})-(\ref{Evolution-ODE-Hybrid-blend}) differs from the initial-boundary-value problem (\ref{Equilibrium-PDE-Hybrid-0})-(\ref{Evolution-ODE-Hybrid-0}) on two counts: ($i$) the material contants are $\bfX$-dependent and ($ii$) the fields $\bfy(\bfX,t)$ and $q(\bfX,t)$ are subject to a periodicity constraint. In order to handle the first of these two differences when using Abaqus, it suffices to make use of a \emph{conforming} FE mesh, that is, a mesh wherein every element is fully comprised of either the rubber phase $r=1$ or the rubber phase $r=2$ as dictated by the characteristic functions (\ref{theta-blend}). To handle the second difference, the FE mesh has to be not only conforming but also \emph{periodic} so that the degrees of freedom in (\ref{FE-Representation-y}) and (\ref{FE-Representation-q}) on the boundary of the unit cell can be constrained appropriately in a manner consistent with the way Abaqus enforces linear multipoint constraints.
\end{remark}

\subsection{Sample results}

\begin{figure}[H]
\centering
\includegraphics[width=0.49\textwidth]{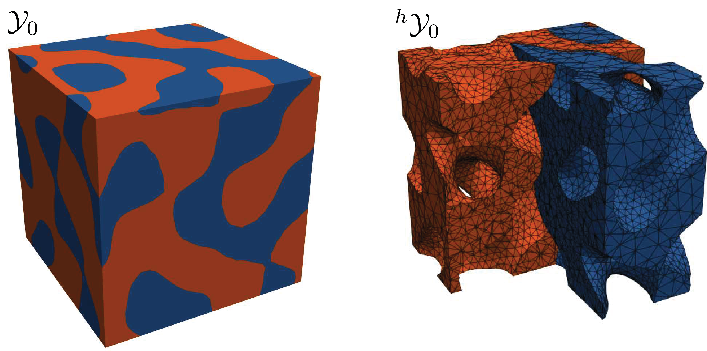}
\caption{Unit cell $\mathcal{Y}_0$ of a bicontinuous rubber blend with equal (50/50) volume fraction of each rubber phase and its FE discretization ${}^h\mathcal{Y}_0$ with about 75,000 simplicial elements. The mesh is clipped in order to better illustrate the bicontinuous character of the microstructure.}
\label{Fig5}
\end{figure}

In a recent contribution \cite{SLPLP24}, we have developed a novel meshing scheme that is capable of generating the type of conforming and periodic meshes --- starting from a voxelized representation of the microstructures of blends --- required to solve (\ref{Equilibrium-PDE-Hybrid-blend})-(\ref{Evolution-ODE-Hybrid-blend}) in Abaqus. Figure \ref{Fig5} shows a mesh generated by that scheme for a representative unit cell describing the rubber blend that we shall study here.
%

%
%%%%%%%%%%%%%%%%%%%%%%%%%%%%%%%%%%%%%%%%%%%%%%%%%%%%%%%%%%%%%%%%%%%%%%%%%%%%%%
\begin{figure}[t!]
  \subfigure[]{
   \label{fig:2a}
   \begin{minipage}[]{0.5\textwidth}
   \centering \includegraphics[width=2.55in]{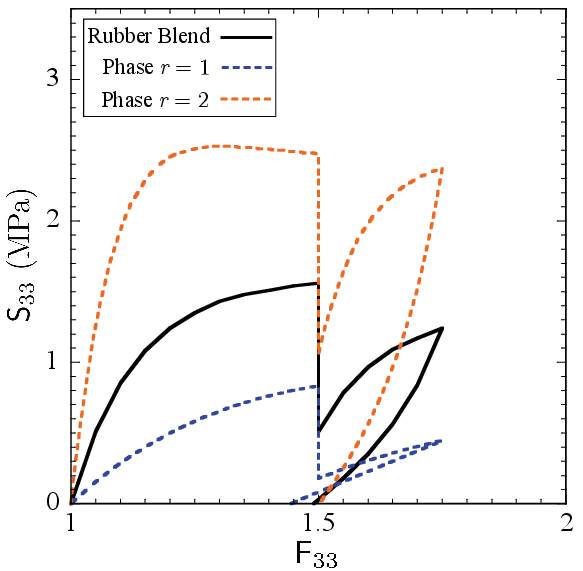}
   \vspace{0.2cm}
   \end{minipage}}
  \subfigure[]{
   \label{fig:2b}
   \begin{minipage}[]{0.5\textwidth}
   \centering \includegraphics[width=2.55in]{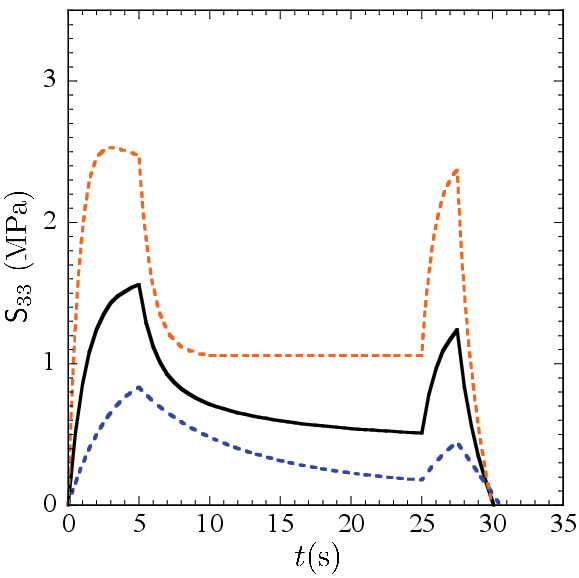}
   \vspace{0.2cm}
   \end{minipage}}
   \caption{Macroscopic response of the bicontinuous rubber blend with the microstructure shown in Fig. \ref{Fig5} --- made of two Gaussian rubber phases with constant viscosity and the material constants listed in Table \ref{Table2} --- under the uniaxial tension loading/relaxation/loading/unloading (\ref{S33-blend})--(\ref{F33-blend}). The results show the macroscopic stress $\mathsf{S}_{33}$: (a) as a function of the prescribed macroscopic deformation $\mathsf{F}_{33}$ and (b) as a function of time $t$. For direct comparison, the corresponding responses of the underlying rubber phases are also plotted.}\label{Fig6}
\end{figure}
%%%%%%%%%%%%%%%%%%%%%%%%%%%%%%%%%%%%%%%%%%%%%%%%%%%%%%%%%%%%%%%%%%%%%%%%%%%%%%
%

At this stage, we are in a position to deploy the proposed UMAT subroutine to solve the unit-cell problem (\ref{Equilibrium-PDE-Hybrid-blend})-(\ref{Evolution-ODE-Hybrid-blend}) in Abaqus and, in turn, compute the resulting history of the macroscopic stress (\ref{S-avg-Y}) in terms of the history of the macroscopic deformation gradient tensor (\ref{F-avg-Y}).

\begin{table}[h!]\centering
\caption{Values of the material constants for the two rubber phases $r=1$ and $r=2$ in the blend.}
\begin{tabular}{llll}
\toprule
$\mu^{(1)}=0.1$ MPa & & & $\kappa^{(1)}=10^4\mu^{(1)}=10^3$ MPa \\
$m^{(1)}=1$ MPa & & & $\eta^{(1)}=10$ MPa s \\
\midrule
\midrule
$\mu^{(2)}=1$ MPa & & & $\kappa^{(2)}=10^4\mu^{(2)}=10^4$ MPa \\
$m^{(2)}=10$ MPa & & & $\eta^{(2)}=1$ MPa s \\
\bottomrule
\end{tabular} \label{Table2}
\end{table}

Figure \ref{Fig6} presents results for the macroscopic response of the blend with the microstructure characterized by the unit cell presented in Fig. \ref{Fig5}, made of Gaussian rubbers with constant viscosity and the material constants listed in Table \ref{Table2}, that is subjected to uniaxial tension loading/relaxation/loading/unloading
\begin{equation}\label{S33-blend}
\boldsymbol{\mathsf{S}}(t)=\mathsf{S}_{33}(t)\bfe_3\otimes\bfe_3
\end{equation}
with
\begin{equation}\label{F33-blend}
\mathsf{F}_{33}(t)=\left\{\begin{array}{ll}
\hspace{-0.15cm}1+\dot{\lambda}_0t, & 0\leq t\leq 5\,{\rm s}\\
\hspace{-0.15cm}1.5, & 5\,{\rm s}< t\leq 25\,{\rm s} \\
\hspace{-0.15cm}\dot{\lambda}_0 t-1, & 25\,{\rm s}< t\leq 27.5\,{\rm s}\\
\hspace{-0.15cm}4.5-\dot{\lambda}_0 t, & 27.5\,{\rm s}<t\leq 35\,{\rm s}
\end{array}\right., \quad\dot{\lambda}_0=10^{-1}\,{\rm s}^{-1},
\end{equation}
where, again, $\{\bfe_1,\bfe_2,\bfe_3\}$ stands for the Cartesian laboratory frame of reference. Specifically, Fig. \ref{Fig6}(a) shows the macroscopic stress $\mathsf{S}_{33}$ as a function of the prescribed macroscopic deformation $\mathsf{F}_{33}$, while Fig. \ref{Fig6}(b) shows $\mathsf{S}_{33}$ as a function of time $t$. For direct comparison, the figures include the corresponding responses of the two rubber phases $r=1$ and $r=2$ that make up the blend.

The type of results presented in Fig. \ref{Fig6} clearly demonstrate the capabilities of the proposed UMAT subroutine to solve complex mechanics problems in finite viscoelastostatics.

The results \emph{per se} also suggest that the effective viscosity of a rubber blend is functionally much richer than the viscosities of the underlying rubber phases, a behavior that is consistent with recent results for suspensions of rigid inclusions and vacuous bubbles in viscoelastic rubber \cite{GSKLP21,SGLP23}. The study of this phenomena will be the subject of future work.

\section*{Acknowledgements}

\noindent Support for this work by the National Science Foundation through the Grant DMS--2308169 is gratefully acknowledged.

\bibliographystyle{elsarticle-num}
\bibliography{refs}

\end{document}